\documentclass[12pt]{article}
\usepackage[dvips]{epsfig}
\usepackage{graphicx}
\usepackage{amsthm,amsmath,amssymb}
\usepackage{mathtools}
\usepackage{latexsym}
\usepackage{graphics}
\usepackage{url}
\usepackage{placeins}
\usepackage{color}
\usepackage{bm}
\usepackage{comment}
\usepackage{statex}
\usepackage{pgfplots}

\usepackage[
            bookmarksnumbered=true,
            bookmarksopen=true,
            colorlinks=true,
            citecolor=green,
            linkcolor=green,
            anchorcolor=red,
            urlcolor=blue
            ]{hyperref}

\definecolor{blue}{rgb}{0,0,0.9}
\definecolor{red}{rgb}{0.9,0,0}
\definecolor{green}{rgb}{0,0.50,0.10}
\definecolor{violet}{rgb}{0.5804,0.0000,0.8275}

\def\@themcountersep{}

\makeatletter
\newcommand{\labeltext}[2]{%
  \@bsphack
  \csname phantomsection\endcsname 
  \def\@currentlabel{#1}{\label{#2}}%
  \@esphack
}
\makeatother

\newtheorem{THEO}{Theorem}[section]
\newtheorem{ALGo}[THEO]{Algorithm}
\newtheorem{ASSUMPT}[THEO]{Assumption}
\newtheorem{CONJ}[THEO]{Conjecture}
\newtheorem{COND}[THEO]{Condition}
\newtheorem{CORO}[THEO]{Corollary}
\newtheorem{DEFI}[THEO]{Definition}
\newtheorem{EXAMP}[THEO]{Example}
\newtheorem{INSTANCE}[THEO]{Instance}
\newtheorem{FACT}[THEO]{Fact}
\newtheorem{HYPO}[THEO]{Hypothesis}
\newtheorem{LEMM}[THEO]{Lemma}
\newtheorem{PROB}[THEO]{Problem}
\newtheorem{PROP}[THEO]{Proposition}
\newtheorem{REMA}[THEO]{Remark}
\newcommand{\theo}{\begin{THEO}}
\newcommand{\algo}{\begin{ALGo} \rm}
\newcommand{\assumpt}{\begin{ASSUMPT} \rm}
\newcommand{\cond}{\begin{COND}}
\newcommand{\conj}{\begin{CONJ}}
\newcommand{\coro}{\begin{CORO}}
\newcommand{\defi}{\begin{DEFI} \rm}
\newcommand{\examp}{\begin{EXAMP} \rm}
\newcommand{\instan}{\begin{INSTANCE} \rm}
\newcommand{\fact}{\begin{FACT}}
\newcommand{\hypo}{\begin{HYPO} \rm}
\newcommand{\lemm}{\begin{LEMM}}
\newcommand{\prob}{\begin{PROB} \rm}
\newcommand{\prop}{\begin{PROP}}
\newcommand{\rema}{\begin{REMA} \rm}

\newcommand{\etheo}{\end{THEO}}
\newcommand{\ealgo}{\end{ALGo}}
\newcommand{\eassumpt}{\end{ASSUMPT}}
\newcommand{\econd}{\end{COND}}
\newcommand{\econj}{\end{CONJ}}
\newcommand{\ecoro}{\end{CORO}}
\newcommand{\edefi}{\end{DEFI}}
\newcommand{\eexamp}{\end{EXAMP}}
\newcommand{\einstan}{\end{INSTANCE}}
\newcommand{\efact}{\end{FACT}}
\newcommand{\ehypo}{\end{HYPO}}
\newcommand{\elemm}{\end{LEMM}}
\newcommand{\eprob}{\end{PROB}}
\newcommand{\eprop}{\end{PROP}}
\newcommand{\erema}{\end{REMA}}

\def\0{\mbox{\bf 0}}
\def\1{\mbox{\bf 1}}
\def\2{\mbox{\bf 2}}
\def\3{\mbox{\bf 3}}
\def\4{\mbox{\bf 4}}
\def\5{\mbox{\bf 5}}
\def\6{\mbox{\bf 6}}
\def\7{\mbox{\bf 7}}
\def\8{\mbox{\bf 8}}
\def\9{\mbox{\bf 9}}
\def\a{\mbox{\boldmath $a$}}
\def\b{\mbox{\boldmath $b$}}
\def\cc{\mbox{\boldmath $c$}}
\def\d{\mbox{\boldmath $d$}}
\def\e{\mbox{\boldmath $e$}}

\def\u{\mbox{\boldmath $u$}}

\def\x{\mbox{\boldmath $x$}}
\def\y{\mbox{\boldmath $y$}}
\def\z{\mbox{\boldmath $z$}}
\def\A{\mbox{\boldmath $A$}}
\def\B{\mbox{\boldmath $B$}}
\def\C{\mbox{\boldmath $C$}}

\def\F{\mbox{\boldmath $F$}}

\def\H{\mbox{\boldmath $H$}}

\def\O{\mbox{\boldmath $O$}}

\def\Q{\mbox{\boldmath $Q$}}

\def\U{\mbox{\boldmath $U$}}

\def\X{\mbox{\boldmath $X$}}
\def\Y{\mbox{\boldmath $Y$}}

\def\AC{\mbox{$\cal A$}}
\def\BC{\mbox{$\cal B$}}
\def\CC{\mbox{$\cal C$}}
\def\DC{\mbox{$\cal D$}}

\def\QC{\mbox{$\cal Q$}}

\def\inprod#1#2{\langle#1, \, #2\rangle}
\def\Inprod#1#2{{\Large \langle}#1, \, #2{\Large \rangle}}

\def\s0{\mbox{\scriptsize \boldmath $0$}}

\def\coneC{\mathbb{C}}

\def\coneCPP{\mathbb{CP}}
\def\coneDNN{\mathbb{DN}}
\def\coneF{\mathbb{F}}
\def\coneG{\mathbb{G}}
\def\coneN{\mathbb{N}}

\def\coneC{\mathbb{C}}
\def\coneJ{\mathbb{J}}
\def\coneK{\mathbb{K}}

\def\Real{\mathbb{R}}

\def\spaceH{\mathbb{H}}
\def\SymMat{\mathbb{S}}

\def\bGamma{\mbox{\boldmath $\Gamma$}}

\begin{document}

\title{Extending Exact Convex Relaxations of Quadratically Constrained Quadratic Programs}

\author{
\normalsize
Masakazu Kojima\thanks{Department of Data Science for Business Innovation, Chuo University, Tokyo, Japan ({\tt kojima@is.titech.ac.jp}).} \and \normalsize
Sunyoung Kim\thanks{Department of Mathematics, Ewha W. University, 
Seoul, 
Korea 
			({\tt skim@ewha.ac.kr}). 
			 The research was supported  by   NRF 2021-R1A2C1003810.} \and \normalsize
Naohiko Arima\thanks{
	({\tt nao$\_$arima@me.com}).} 
}

\date{\today}

\maketitle

\begin{abstract}

\noindent
A convex 
relaxation of a quadratically constrained quadratic program 
(QCQP) is called exact if it has a rank-$1$ optimal solution that corresponds to 
an optimal solution of the  
QCQP. Given a QCQP whose convex relaxation is exact, 
this paper investigates the incorporation of 
additional quadratic inequality constraints 
under a non-intersecting quadratic constraint condition 
while maintaining the exactness of the convex relaxation of the resulting QCQP. 
Specifically, we extend existing exact semidefinite programming  relaxation, completely positive 
programming  relaxation and doubly nonnegative programming  relaxation of 
various classes of QCQPs in a unified manner. 
Illustrative examples are included to demonstrate the applicability of the established result.
\end{abstract}


\vspace{0.5cm}


\section{Introduction}

\label{section:Introduction}

The quadratically constrained quadratic program (QCQP) seeks to minimize a quadratic function in 
real variables over the feasible region defined by quadratic inequalities. 
The problem is known to be NP-hard \cite{MURTY87}. 
Various convex relaxations  have been extensively studied 
both as theoretical tools and as numerical solution methods for the QCQP. 
Notable examples include the semidefinite programming (SDP) relaxation 
\cite{FUJIE1997,SHOR1987,ZHANG2000}, the doubly nonnegative programming 
(DNN) relaxation \cite{KIM2013,YOSHISE2010}, 
and the completely programming (CPP) relaxation 
\cite{BOMZE2000,BURER2009}. 
 In general, the optimal value $\varphi$ of the QCQP is bounded  by the optimal 
 value $\psi$ of its convex relaxation from below. 
 We say that the QCQP and the convex relaxation are {\em equivalent} if $\psi = \varphi$.
  The convex relaxation is {\em solvable} if it has an optimal solution, 
 and {\em exact} if  it has a rank-$1$ optimal solution that 
 corresponds to an optimal solution of the QCQP with $\psi = \varphi$. 
 
\medskip

This paper investigates extending a given conic-form (or geometric-form) QCQP 
(whose precise description is given in Section~\ref{section:conicQCQP}) 
with an exact convex relaxation by adding quadratic inequality constraints.  
Specifically, the main theorem, Theorem~\ref{theorem:main1}, provides a set of sufficient conditions 
(C1) - (C4) for the convex relaxation of the extended QCQP to remain exact. 
Condition (C1) requires that the convex relaxation of the given QCQP is exact if it is solvable, 
and condition (C4) that the convex relaxation of the resulting QCQP has an optimal solution 
satisfying the KKT (Karush-Kuhn-Tucker) stationary condition. 
Conditions (C2) and (C3) characterize the quadratic inequality constraints to be added.

\medskip

Condition (C1), together with the conic-form for QCQP representation, 
offers a  general and flexible framework that covers  a broad class of 
QCQPs with exact convex relaxations. 
In particular, we  consider the following classes as examples. 
\begin{description} \vspace{-2mm}
\item{(a)} 
A class of QCQPs characterized by  non-intersecting quadratic constraint (NIQC) conditions  
\cite{ARIMA2023,ARIMA2024,BECK2006,BURER2015,CONN2000,JEYAKUMAR2014,JOYCE2024,PONG2014,YANG2018}. 
See Section~\ref{section:NIQC}. 
\vspace{-3mm}
\item{(b) }  
A class of QCQPs characterized by the 
ROG (rank-one-generated) cone condition  \cite{ARGUE2023,ARIMA2023,ARIMA2024,KIM2020}. 
See Section~\ref{section:ROG}. 
\vspace{-3mm}
\item{(c) }  
A class of convex QCQPs. 
See Section~\ref{section:convex}.
\vspace{-2mm}
\item{(d) }  
A class of QCQPs characterized by the  sign pattern condition 
 \cite{AZUMA2022,AZUMA2023,KIM2003,SOJOUDI2014}. 
 See Section~\ref{section:signPattern}.
 \vspace{-3mm}
\item{(e) } 
A continuous quadratic submodular minimization problem 
from \cite{BURER2025}.  
See Section~\ref{section:submodular}.
\vspace{-7mm}
\item{(f) } 
A class of combinatorial QCQPs  \cite{BURER2009,KIM2020}.
See Section~\ref{section:binary}.
\vspace{-3mm}
\item{(g) }  A class of QCQPs with the simplex constraint \cite{GOKMEN2022}. See 
Section~\ref{section:standardQOP}.
\vspace{-2mm}
\end{description}
In classes (a), (b) and (f), the objective function may be an arbitrary quadratic function, 
whereas in the other classes, certain conditions are imposed on both objective and constraint quadratic functions. 
The SDP relaxation is used in (a) through (e), the CPP  relaxation in (f) and the DNN relaxations in (g).
It was shown in \cite[Theorems 1.5, 1.6 and 1.7]{ARIMA2024} that the common NIQC conditions are 
special cases of the ROG cone condition. See Figure~1, where~\eqref{eq:condArimaMain1} 
represents the homogenized NIQC condition proposed by the authors \cite{ARIMA2024} for the conic-form QCQP, 
and \eqref{eq:NICQ1} and~\eqref{eq:NICQ2} denote the non-homogenized NIQC conditions 
\cite{BECK2006,BURER2015,CONN2000,JEYAKUMAR2014,JOYCE2024,PONG2014,YANG2018} 
for the standard-form QCQP (whose precise description is given in 
 Section~\ref{section:standardQCQP}). 
Accordingly, we regard (a) as a subclass of (b). 
The details of the homogenized and non-homogenized 
NIQC conditions, and their relationship are presented in Section~\ref{section:NIQC}. 
Conditions (C2) and (C3) assumed in Theorem~\ref{theorem:main1} to 
characterize the quadratic inequality constraints to be added are variants of the homogenized NIQC 
condition 
~\eqref{eq:condArimaMain1}.

\medskip

Both non-homogenized and homogenized NIQC conditions have been  extensively investigated 
for exact SDP relaxations of QCQPs, as mentioned in (a) above. 
In particular, we refer to the works of  Yang, Anstreiher and Burer  \cite{YANG2018}, 
and Joyse and Yang \cite{JOYCE2024}.  
The aim of these two papers is in line with our main Theorem~\ref{theorem:main1} 
in the sense that they  also address an extension of a given standard-form QCQP 
with an exact SDP relaxation under certain assumptions including non-homogenized 
NIQC conditions on quadratic inequality constraints to be added. 
Since our study covers exact convex relaxations (including SDP, DNN and 
CPP relaxations) of conic-form QCQPs, the QCQPs addressed here are more general. 
Moreover, even if we restrict our main Theorem~\ref{theorem:main1} to exact SDP relaxations of 
standard-form QCQPs, it provides significant improvements over their main results,
\cite[Theorem~1]{YANG2018} and \cite[Corollary~2]{JOYCE2024}. 
The quadratic inequality constraint to be added 
in \cite[Theorem~1]{YANG2018} is restricted to `non-intersecting ellipsoidal hollows'. 
Specifically, if each inequality is expressed as $q(\u) \geq 0$ for some quadratic function $q$, then 
the set $\{ \u : q(\u) \leq 0 \}$ represents an ellipsoid. This assumption is more restrictive and 
less general 
than the NIQC condition assumed in condition (C3), although \cite[Theorem 1]{YANG2018} remains applicable to SDP relaxations of classes (a) - (e) of QCQPs. 

\medskip

Corollary 2 of Joyse and Yang \cite{JOYCE2024} can also be readily applied to extend 
a given standard-form QCQP with the exact SDP relaxation 
as  will be shown in Theorem~\ref{theorem:JOYCE2024a}, 
 although this was not  mentioned explicitly in their paper.  In this case, the exactness of the SDP relaxation of the given QCQP is required for
any choice of its quadratic objective function. In contrast, our condition (C1) (and \cite[Theorem 1]{YANG2018}) 
 requires only that the SDP relaxation of a given QCQP with a fixed quadratic objective function 
be exact when it is solvable. 
Therefore, our requirement on the given QCQP 
is considerably weaker than theirs, a difference that is critical in applications.  
For example, their Corollary 2 cannot be applied to class (c) of convex QCQPs and 
(d) of QCQPs characterized by the sign pattern condition. 
They also assumed a variant of non-homogenized NIQC condition~\eqref{eq:NICQ2} for 
quadratic inequality constranits to be added, which 
is slightly different from condition (C3), a variant of homogenized NIQC 
condition~\eqref{eq:condArimaMain1}. See Figure~1.

\medskip

 The main contributions of this paper are summarized as follows:
 \begin{itemize} \vspace{-2mm}
  \item We present  a unified framework that extends a given 
  QCQP with an exact convex relaxation by adding quadratic inequalities satisfying 
 non-homogenized NIQC conditions. 
 Since condition (C1), imposed on a conic-form representation of the given QCQP, 
 is quite general and flexible, our framework can be applied to a wide class of QCQPs 
 (represented not only in the conic-form but also in the standard-form) with exact convex relaxations. 
\vspace{-2mm}
\item Those classes include classes (b) - (e) of QCQPs with exact SDP relaxations, class (f) of QCQPs with 
 exact CPP relaxations and class (g) of QCQPs with exact DNN relaxations, which have been investigated 
largely independently in the existing literature.
 \vspace{-2mm}
  \item Furthermore, our unified framework can be applied to any convex relaxation over 
 more general closed convex cone $\coneK$  contained in the positive semidefinite cone,
 if it is  known to result in exact convex relaxations.
  \vspace{-2mm}
 \end{itemize}
 
\subsection*{Outline of the paper}

After describing some notation and symbols used throughout the paper, we introduce 
two distinct representations of QCQPs: a conic-form QCQP in Section~\ref{section:conicQCQP} 
and a standard-form QCQP in Section~\ref{section:standardQCQP}. 
The standard-form QCQP is a special case of the conic-form QCQP, and 
the convex relaxations including SDP, CPP and DNN relaxations are described for the conic-form QCQP in Section~\ref{section:conicQCQP}. 
In Section~\ref{section:NIQC}, various NIQC conditions are stated and their relationship are shown. 
Our main theorem, Theorem~\ref{theorem:main1}, which includes conditions (C1) - (C4),  
is given in Section~\ref{section:MainTheorem}. 
Also some 
 existing results related to
 Theorem~\ref{theorem:main1}, 
including Theorem~\ref{theorem:JOYCE2024a} 
obtained from \cite[Corollaries 2]{JOYCE2024} mentioned above as a special case, 
are presented. 
Theorem~\ref{theorem:main1} is applied 
to exact SDP relaxations of QCQPs from the classes (b), (c), (d) and (e) in Section~\ref{section:SDP}, 
and to exact CPP and DNN relaxation of QCQPs from the classes (f) and (g) in Section~\ref{section:CPPandDNN}. 
We conclude the paper in Section~\ref{section:Conclusion}.

\section{Preliminaries}

We use the following notation: 
\begin{eqnarray*}
& & \Real^n : \mbox{the $n$-dimensional linear space of column vectors $\x = (x_1,\ldots,x_n)$}, \\[3pt]
& &  \Real^n_+ = \{\x \in \Real^n : \x \geq \0\} \ \mbox{(the nonnegative orthant of $\Real^n$)}, \\[3pt] 
& & \SymMat^n : \mbox{the linear space of $n \times n$ symmetric matrices}, \\
& & \inprod{\A}{\X}  =  \sum_{i=1}^n\sum_{j=1}^n A_{ij}X_{ij} \ \mbox{(the inner product of $\A, \ \X \in 
                                 \SymMat^n$)} \ 
\mbox{for every } \A,\X \in \SymMat^n, \\[3pt]
& &  \SymMat^n_+ : \mbox{the convex cone of $n \times n$ positive semidefinite symmetric matrices},
\\[3pt]
& &  \bGamma^n = \{ \x\x^T : \x \in \Real^n\}, \ \mbox{where $\x^T$ denotes the transposed row vector of $\x \in \Real^n$}, \\[3pt]
& & \mbox{co}S \ \mbox{and }  \overline{\mbox{co}}S : 
\mbox{the convex hull of $S \subseteq \SymMat^n$ and its closure, respectively}.
\end{eqnarray*}
 For every $\A \in \SymMat^n$, 
the quadratic form $\x^T\A\x$ in $\x \in \Real^n$ is frequently expressed as $\inprod{\A}{\x\x^T}$. 

\subsection{A conic-form QCQP and its convex relaxation}

\label{section:conicQCQP}

Let $\Q \in \SymMat^n$ and  $\O \not=\H \in \SymMat^n$. 
For every closed cone $\coneC \subset \SymMat^n_+$,  
COP$(\coneC,\Q,\H)$ denotes 
the problem of minimizing $\inprod{\Q}{\X}$ subject to 
$\X \in \coneC$ and $\inprod{\H}{\X} = 1$, {\it i.e.}, 
\begin{eqnarray*}
\eta(\coneC,\Q,\H) & = & \inf \left\{ \inprod{\Q}{\X} : \X \in \coneC, \ \inprod{\H}{\X} = 1 \right\}.
\end{eqnarray*}
Here $\coneC \subseteq \SymMat^n_+$  is a cone if $\lambda\X \in \coneC$ for every $\lambda \geq 0$  and $\X \in \coneC$.
We note that every feasible solution $\X$ of COP$(\coneC,\Q,\H)$ is nonzero since $\H \not= \O$.
If COP$(\coneC,\Q,\H)$ is infeasible, we assume that $\eta(\coneC,\Q,\H) = +\infty$. 
We say that COP$(\coneC,\Q,\H)$ is solvable if it has an optimal solution.
For every closed convex cone $\coneG \subseteq \SymMat^n_+$, 
{\em the conic-form (geometric-form)} QCQP and its convex relaxation 
are described as COP$(\bGamma^n\cap\coneG,\Q,\H)$ 
and COP$(\coneG,\Q,\H)$, 
respectively. 
Obviously, $\bGamma^n\cap\coneG = \{\X \in \coneG: \mbox{rank}\X \leq 1 \}$ and 
$\eta(\coneG,\Q,\H) \leq \eta(\bGamma^n\cap\coneG,\Q,\H)$. 
We say that COP$(\coneG,\Q,\H)$ and COP$(\bGamma^n\cap\coneG,\Q,\H)$ are {\em equivalent} 
if $\eta(\coneG,\Q,\H) = \eta(\bGamma^n\cap\coneG,\Q,\H)$, and that COP$(\coneG,\Q,\H)$ 
is {\em an exact convex relaxation} (or simply {\em exact}) if 
it has an optimal solution $\X\in\bGamma^n$, which is an optimal solution of 
COP$(\bGamma^n\cap\coneG,\Q,\H)$.

\medskip

We now focus on three main types of convex relaxation, the SDP relaxation, 
the CPP relaxation,  and  the DNN relaxation. 
The cones associated with these relaxations are $\SymMat^n_+$, %
\begin{eqnarray*} 
\coneCPP^n & = & \mbox{co}\{\x\x^T : \x \in \Real^n_+\} \
 \mbox{(the CPP cone)}, \\[3pt]
 \coneDNN^n & = &  \SymMat^n_+ \cap \coneN^n 
\ \mbox{(the DNN cone), where } \\
&\  &  \coneN^n = \{ \X \in \SymMat^n : X_{ij} \geq 0 \ (1 \leq i, \ j\leq n)\}, 
\end{eqnarray*}
respectively. 
When $\coneK \in \{\coneDNN^n,\coneCPP^n\}$, $\x\x^T \in \bGamma^n\cap\coneK$ if and only if 
$\x \in \Real^n_+\cup(-\Real^n_+)$. Since $\x\x^T = (-\x)(-\x)^T$ for every $\x \in \Real^n$, 
we may restrict $\x \in \Real^n_+$ in COP$(\bGamma^n\cap\coneG,\Q,\H)$ with $\coneG \subseteq 
\coneK \in \{\coneDNN^n,\coneCPP^n\}$. Thus, we can write COP$(\bGamma^n\cap\coneG,\Q,\H)$ as 
QCQP 
\begin{eqnarray}
\eta(\bGamma^n\cap\coneG,\Q,\H) = \inf\left\{\inprod{\Q}{\x\x^T}: 
\begin{array}{l}
\x \in \Real^n \ \mbox{if } \coneG \subseteq \coneK = \SymMat^n_+, \\
\x \in \Real^n_+ \ \mbox{if } \coneG \subseteq \coneK \in  \{\coneDNN^n,\coneCPP^n\}, \\
\x\x^T \in \coneG, \ \inprod{\H}{\x\x^T} = 1
\end{array}
\right\}. \label{eq:QCQPG}
\end{eqnarray}

\medskip

We often represent $\coneG$ 
using a closed convex cone 
$\coneK \subseteq \SymMat^n_+$ and linear matrix inequalities in $\X \in \SymMat^n_+$. 
For a simple description of  a closed convex cone determined by linear matrix inequalities and equalities, 
we define 
\begin{eqnarray*}
& & \coneJ_-(\A), \ \coneJ_0(\A) \ \mbox{or } \coneJ_+(\A)  \equiv  
\left\{ \X \in \SymMat^n_+ : \inprod{\A}{\X} \leq, \  = \ \mbox{or } \geq 0, \ \mbox{respectively} 
\right\} ,
\\[3pt] 
& & \coneJ_+(\AC)  \equiv  \bigcap_{\A \in \AC} \coneJ_+(\A)
 = \left\{ \X \in \SymMat^n_+ : \inprod{\A}{\X} \geq 0 \ 
(\A \in \AC) \right\} 
\end{eqnarray*}
for  every $\A \in \SymMat^n$ and $\AC \subseteq \SymMat^n$. 
Now, letting $\coneG = \coneK \cap \coneJ_+(\AC)$ for some finite $\AC \subset \SymMat^n$, we rewrite 
QCQP~\eqref{eq:QCQPG} as 
\begin{eqnarray}
\eta(\bGamma^n\cap\coneK\cap\coneJ_+(\AC),\Q,\H) 
& = & \inf\left\{\inprod{\Q}{\x\x^T}: 
\begin{array}{l}
\x \in \Real^n \ \mbox{if } \coneK = \SymMat^n_+, \\
\x \in \Real^n_+ \ \mbox{if } \coneK \in  \{\coneDNN^n,\coneCPP^n\}, \\
\inprod{\A}{\x\x^T} \geq 0  \ (\A \in \AC), \\ 
\inprod{\H}{\x\x^T} = 1
\end{array}
\right\}\nonumber  \\[3pt]
& = & \inf\left\{\inprod{\Q}{\x\x^T}: 
\begin{array}{l}
\x \in \Real^n \ \mbox{if } \coneK = \SymMat^n_+, \\
\x \in \Real^n_+ \ \mbox{if } \coneK \in  \{\coneDNN^n,\coneCPP^n\}, \\
\x\x^T \in \coneJ_+(\AC), \ \inprod{\H}{\x\x^T} = 1
\end{array}
\right\}, \label{eq:QCQPLMI}
\end{eqnarray}
and its convex relaxation COP$(\coneK\cap\coneJ_+(\AC),\Q,\H))$ as 
\begin{eqnarray*}
\eta(\coneK \cap \coneJ_+(\AC),\Q,\H) & = & \inf\left\{ \inprod{\Q}{\X} : \X \in \coneK\cap\coneJ_+(\AC), \
\inprod{\H}{\X} = 1 
\right\}. 
\end{eqnarray*}

\subsection{A standard-form QCQP}

\label{section:standardQCQP} 

We present a standard-form QCQP for minimizing a quadratic objective function subject to 
quadratic inequalities in a finite number of (nonnegative) real variables, which is more commonly used than the conic-form QCQP~\eqref{eq:QCQPLMI}. The standard-form QCQP is obtained from the conic-form QCQP~\eqref{eq:QCQPLMI} by fixing $\H = \H^1 \equiv \mbox{diag}(0,0,\ldots,1) \in \SymMat^n_+$ 
(the $n \times n$ diagonal matrix with the diagonal elements $0,0,\ldots,1$). 
Then, the constraint $\inprod{\H}{\x\x^T} = 1$ can be replaced by 
$x_n \in \{-1,1\}$. Since $\x\x^T = (-\x)(-\x)^T$ for every $\x \in \Real^n$, we may fix $x_n=1$. 
We denote each $\x \in \Real^n$ with $x_n = 1$ as 
$\x = {\scriptsize \begin{pmatrix} \u\\ 1\end{pmatrix}}$ 
with $\u \in \Real^{n-1}$. Then, in QCQP~\eqref{eq:QCQPLMI} with $\H = \H^1$,  
the quadratic form $\inprod{\A}{\x\x^T}$ of 
$\x$ with $x_n = 1$ 
can be replaced by 
 a quadratic function in $\u \in \Real^{n-1}$; 
\begin{eqnarray*}
q(\u,\A) & = &
\inprod{\A}{{\scriptsize \begin{pmatrix}\u\\1\end{pmatrix}}{\scriptsize \begin{pmatrix}\u\\1\end{pmatrix}}^T}
 \ \mbox{for every } \u \in \Real^{n-1}. 
\end{eqnarray*}
Define 
subsets $\A_\leq, \ \A_=, \ \A_\geq, \ \AC_\geq$ of $\Real^{n-1}$ by 
\begin{eqnarray*}
& & \A_\leq, \ \A_= \ \mbox{or } \A_\geq = \left\{ \u \in \Real^{n-1} : q(\u,\A) \leq, \ = \ \mbox{or } \geq 0, 
\ \mbox{respectively}\right\}, \\
& & \AC_\geq =  \bigcap_{\scriptsize \A \in \AC} \A_\geq = \left\{ \u \in \Real^{n-1} : q(\u,\A) \geq 0 \ 
(\A \in \AC) \right\}. 
\end{eqnarray*}
As a result, QCQP~\eqref{eq:QCQPLMI} with $\H = \H^1$ can be rewritten as 
{\em a standard-form} QCQP 
\begin{eqnarray}
\eta(\bGamma^n\cap\coneK\cap\coneJ_+(\AC),\Q,\H^1) 
& = & \inf \left\{ q(\u,\Q) : 
\begin{array}{l}
\u \in \Real^{n-1} \ \mbox{if } \coneK = \SymMat^n_+, \\
\u \in \Real^{n-1}_+ \ \mbox{if } \coneK \in  \{\coneDNN^n,\coneCPP^n\}, \\
q(\u,\A) \geq 0 \ (\A \in \AC)
\end{array}
\right\} \nonumber \\[3pt]
& = & \inf \left\{ q(\u,\Q) : 
\begin{array}{l}
\u \in \Real^{n-1} \ \mbox{if } \coneK = \SymMat^n_+, \\
\u \in \Real^{n-1}_+ \ \mbox{if } \coneK \in  \{\coneDNN^n,\coneCPP^n\}, \\
\u \in \AC_\geq
\end{array}
\right\}. 
 \label{eq:QCQPS}
\end{eqnarray}
It should be noted that QCQP~\eqref{eq:QCQPS} is a special case of QCQP~\eqref{eq:QCQPLMI} 
by fixing $\H = \H^1$. Therefore, they share the common convex relaxation 
COP$(\coneK\cap\coneJ_+(\AC),\Q,\H^1)$.

\subsection{Non-intersecting quadratic constraint (NIQC) conditions }

\label{section:NIQC}

Throughout this section, we assume that $\coneK = \SymMat^n_+$. 
A  NIQC condition for QCQP~\eqref{eq:QCQPLMI} with $\coneK = \SymMat^n_+$ can be stated as 
\begin{eqnarray}
& & \coneJ_0(\A) \subseteq \coneJ_+(\A') \ \mbox{for every distinct } \A, \ \A' \in \AC.  
\label{eq:NIQCLMI0}
\end{eqnarray}
In \cite{ARIMA2023}, this condition was introduced as a sufficient condition for 
$\coneJ_+(\AC)$ to be ROG (rank-one-generated), 
which characterizes the exactness of the SDP relaxation of QCQP~\eqref{eq:QCQPLMI} 
\cite[Theorem~4.1]{ARIMA2023}.

\medskip

By construction, we know that 
\begin{eqnarray*}
& & \A_\leq, \ \A_= \ \mbox{or } \A_\geq   = \left\{ \u \in \Real^{n-1} :  {\scriptsize \begin{pmatrix}\u\\1\end{pmatrix}}{\scriptsize \begin{pmatrix}\u\\1\end{pmatrix}}^T \in \coneJ_-(\A), \  \coneJ_0(\A)  \ 
\mbox{or } \coneJ_+(\A), \ \mbox{respectively} \right\}, \\ 
& &  \AC_\geq  = \left\{ \u \in \Real^{n-1} :  {\scriptsize \begin{pmatrix}\u\\1\end{pmatrix}}{\scriptsize \begin{pmatrix}\u\\1\end{pmatrix}}^T \in \coneJ_+(\AC) \right\}. \label{eq:NIQCS}  
\end{eqnarray*}
Another type of NIQC condition for QCQP~\eqref{eq:QCQPS} with $\coneK = \SymMat^n_+$
\begin{eqnarray}  
\A_= \subseteq \A'_\geq \ \mbox{for every distinct } \A, \ \A' \in \AC \label{eq:NIQCS} 
 \end{eqnarray}
 has been used (implicitly) in much of the 
 existing literature. 
 To distinguish NIQC conditions~\eqref{eq:NIQCLMI0} for QCQP~\eqref{eq:QCQPLMI} 
 and~\eqref{eq:NIQCS} for QCQP~\eqref{eq:QCQPS}, we refer to them  as homogenized and 
 non-homogenized, respectively. 
Condition~\eqref{eq:NIQCS} 
was originally studied for  exact SDP relaxations for  
simple QCQPs, particularly those arising from the generalized trust region subproblem (TRS) 
\cite{BECK2006,BURER2015,CONN2000,JEYAKUMAR2014,PONG2014}. 
As an extension of the generalized TRS, 
the works in \cite{POLYAK1998,YE2003} investigated  
a QCQP of the form $\inf\{q_0(\u) : -1 \leq q_1(\u) \leq 1\}$, 
where $q_0, \ q_1$ are  quadratic functions in $\u \in \Real^{n-1}$. 
This problem can be reformulated as 
a QCQP satisfying NIQC condition. 
Similarly, a quadratic program with non-intersecting ellipsoidal hollows \cite{YANG2018} provides another extension of the generalized TRS. 
Recently, \cite{JOYCE2024} applied the non-homogenized NIQC condition to obtain exact SDP relaxations  
of more general QCQPs, which we  present in Theorem~\ref{theorem:JOYCE2024a}.

\medskip

Generalizations of NIQC conditions~\eqref{eq:NIQCLMI0} 
and~\eqref{eq:NIQCS} 
were proposed in \cite{ARIMA2024} for the cases where $\coneJ_+(\AC)$ is 
contained in a face $\coneF$ of  $\SymMat^n_+$ and $\AC$ possibly contains 
redundant $\A \in \AC$ to represent $\coneF \cap \coneJ_+(\AC)$. For example, this occurs if  
$\AC$ contains a negative semidefinite matrix, and also
in the study CPP relaxation for a class of 
combinatorial QCQPs in Section~\ref{section:binary}.  See also \cite[Section 2]{ARIMA2024}. 

\medskip

In the remainder of this section, we assume that $\coneJ_+(\AC) \subseteq \coneF$ for some face $\coneF$ 
of $\SymMat^n_+$ represented as 
$\coneF = \{\X \in \SymMat^n_+ : \inprod{\F}{\X} = 0 \}$ with $\F \in \SymMat^n_+$. 
Define 
\begin{eqnarray*}
& & L(\coneF) = \{\x \in \Real^n : \x\x^T \in \coneF\} = \{\x \in \Real^n : \inprod{\F}{\x\x^T} = 0\}
= \{\x \in \Real^n :\F\x = \0\}, \\ 
& & L_1(\coneF) = \{ \u \in \Real^{n-1} : {\scriptsize \begin{pmatrix}\u \\ 1\end{pmatrix} }\in L(\coneF) \} 
= \{ \u \in \Real^{n-1} : \F {\scriptsize \begin{pmatrix}\u \\ 1\end{pmatrix} } = \0 \}. 
\end{eqnarray*}
Then, $L(\coneF)$ forms a linear subspace of $\Real^n$ that includes  
$\{\x \in \Real^n: \x\x^T \in \coneJ_+(\AC)\}$, and $L_1(\F)$ an affine subspace of $\Real^{n-1}$ 
that includes  $\AC_\geq$.  It should be noted that $\coneJ_+(\AC) = \coneF\cap\coneJ_+(\AC)$ 
and $\AC_\geq = L_1(\coneF)\cap\AC_\geq$. 

\medskip

\begin{figure}[t!]  
\begin{picture}(400,250)(0,0)
\put(50,200){\framebox(90,20)}
\put(225,190){\framebox(215,40)}
\put(60,205){$\coneJ_+(\AC)$ is ROG}
\put(149,215){Lemma~\ref{lemma:ROG}}
\put(167,205){$\Longrightarrow$}
\put(232,215){$\eta(\coneJ_+(\AC),\Q,\H) = \eta(\bGamma^n\cap\coneJ_+(\AC),\Q,\H)$}
\put(232,198){\hspace{2.5cm}  if $-\infty < \eta(\coneJ_+(\AC),\Q,\H)$}
\put(110,183){$\Uparrow$}
\put(120,181){Theorem~\ref{theorem:ArimaMain1}}
%
\put(37,155){
 \eqref{eq:NIQCLMI0}  \hspace{25mm} \eqref{eq:condArimaMain1} \hspace{45mm} \eqref{eq:ArgueProp1}}
\put(70,155){$\Longrightarrow$}
\put(22,150){\framebox(40,20)}
\put(190,171){ \cite[Section 5.1]{ARIMA2024}}
\put(95,150){\framebox(90,20)}
\put(220,142){$\sharp$}
\put(213,150){$\Longrightarrow$}
\put(213,160){$\Longleftarrow$}
\put(273,150){\framebox(35,20)}
%
\put(118,130){\ $\sharp$ \hspace{6mm} $\sharp$}
\put(95,130){\  $\Uparrow$ \ $\Downarrow$ \hspace{10mm} $\Downarrow$ \ $\Uparrow$}
\put(90,100){\framebox(40,20)}
\put(145,100){\framebox(40,20)}
\put(100,105){\eqref{eq:NICQ1} \hspace{12mm} \eqref{eq:NICQ2}}
\put(08,129){\   \cite[Theorem 1.6]{ARIMA2024}}
\put(178,129){\ \cite[Theorem 1.7]{ARIMA2024}}
\end{picture}
\vspace{-35mm} 
\caption{ Summary of NIQC conditions. $\sharp$ : Some assumptions are necessary to prove the relation, including $\coneF = \coneF_{\min}$ 
(the minimal face of $\SymMat^n_+$ that includes $\coneJ_+(\AC)$); 
see \cite[Section 3, 4 and 5]{ARIMA2024} for details. 
}
\end{figure} 

As a generalization of~\eqref{eq:NIQCLMI0}, we have 
\begin{eqnarray}
& & \coneF  \cap  \coneJ_0(\A)\subseteq  \coneF \cap \coneJ_+(\A')\ \mbox{ \ or \ } 
\coneF \cap \coneJ_+(\A')  \subseteq  \coneF  \cap \coneJ_+(\A) \ \mbox{for every } \A, \ \A' \in \AC
\label{eq:condArimaMain1}
\end{eqnarray}
(\cite[condition (B)]{ARIMA2024}). 
The first inclusion relation  corresponds to the homogenized NIQC condition restricted to $\coneF$. 
The second implies that $\A \in \AC$ is redundant for describing $\coneF\cap\coneJ_+(\AC)$. 
Specifically, $\coneF\cap\coneJ_+(\AC) = \coneF\cap\coneJ_+(\AC\backslash\{\A\})$. However,
$\coneJ_+(\AC)$ could be a proper subset of $\coneJ_+(\AC\backslash\{\A\})$ 
\cite[Section 2]{ARIMA2024}. 
Condition~\eqref{eq:condArimaMain1} ensures that $\coneJ_+(\AC)$ is ROG, and identifies a class of QCQPs whose 
SDP relaxations are exact, as  discussed in Section~\ref{section:ROG}. 
It was shown in \cite[Section 5]{ARIMA2024} that 
condition~\eqref{eq:condArimaMain1} is weaker 
than any of the three variants of the NIQC condition stated below. 

\medskip

The first condition is a non-homogenized NIQC condition assumed in \cite[Theorem 1.5]{ARIMA2024}:  
\begin{eqnarray}
& & 
\left.
\begin{array}{l}
\AC_\geq = L_1(\coneF) \cap\AC_\geq \not= \emptyset \ \mbox{and }\\[3pt]
\emptyset \not=   L_1(\coneF) \cap \A_\leq \subseteq L_1(\coneF) \cap  \A'_\geq 
\ \mbox{or } L_1(\coneF) \cap \A'_\geq   \subseteq   L_1(\coneF) \cap \A_\geq\\
\hspace{80mm} \mbox{for every } \A, \ \A' \in \AC. 
\end{array}
\right\} \label{eq:NICQ1}
\end{eqnarray}
The second one is another non-homogenized NIQC condition assumed in \cite[Theorem 1.6]{ARIMA2024}, which 
is obtained from \cite[Corollary 3]{JOYCE2024}: 
\begin{eqnarray}
\left.
\begin{array}{l}
L_1(\coneF)\cap\A_=  \subseteq L_1(\coneF)\cap\A'_\geq \ \mbox{or } 
L_1(\coneF)\cap\A'_\geq  \subseteq L_1(\coneF)\cap\A_\geq \\
\hspace{80mm} 
\mbox{ for every } \A, \ \A' \in \AC,\\[3pt]
 q(\cdot,\A) : \Real^{n-1} \rightarrow \Real \ \mbox{is not affine on $L_1(\coneF)$} \\
\mbox{({\it i.e.}, the quadratic term of $q(\cdot,\A)|L_1(\coneF)$ is not identically zero)}  \ (\A \in \AC).
\end{array}
\right\} \label{eq:NICQ2}
\end{eqnarray}
The third one is: 
 \begin{eqnarray}
 &  \mbox{for every $\A, \ \A' \in \AC$, there exists $\0 \not= (\alpha,\alpha') \in \Real^2$ such that $\alpha\A+\alpha'\A' \in \SymMat^n_+$} \label{eq:ArgueProp1}
\end{eqnarray}
assumed in  \cite[Proposition 1]{ARGUE2023}. This condition was interpreted as a dual version of condition~\eqref{eq:condArimaMain1} under additional assumptions \cite[Lemma 5.3]{ARIMA2024}. 
Figure 1 summarizes the relationship of the  homogenized and non-homogenized NIQC conditions 
presented above. 

\rema \label{remark:facialReduction}
It is well-known that every face of $\SymMat^n_+$ is isomorphic to $\SymMat^r_+$ for some 
$r \in \{0,\ldots,n\}$ \cite{PATAKI2000,PATAKI2013}. 
Let $\Phi$ be a linear isomorphism from $\coneF \supseteq \coneJ_+(\BC)$ 
onto $\SymMat^r_+$, and 
$\Psi : \SymMat^n \rightarrow \SymMat^r$ the adjoin map with respect to 
$\Phi$.  Then, the pair of COP$(\bGamma^n\cap\coneJ_+(\AC),\Q,\H)$ and its SDP relaxation 
COP$(\coneJ_+(\AC),\Q,\H)$ are equivalently reduced to  the pair of 
COP$(\bGamma^r\cap\coneJ_+(\Psi(\AC)),\Psi(\Q),\Psi(\H))$ and its SDP relaxation 
COP$(\coneJ_+(\Psi(\AC)),\Psi(\Q),\Psi(\H))$, respectively.  
(See \cite{BORWEIN1981,BORWEIN1981a,WAKI2013} 
for numerical methods for the facial  reduction). For the reduced pair, 
\eqref{eq:condArimaMain1} is simplified to 
$
\coneJ_0(\widetilde{\A}) \subseteq  \coneJ_+(\widetilde{\A}') \ \mbox{or } 
 \coneJ_+(\widetilde{\A}')\subseteq  \coneJ_+(\widetilde{\A})
 \ \mbox{for every } 
 \widetilde{\A}, \ \widetilde{\A}' 
 \in \Psi(\AC)$. 
Conditions~\eqref{eq:NICQ1} and ~\eqref{eq:NICQ2} can be simplified accordingly. See Sections 3, 4, and~5 of 
\cite{ARIMA2024} for more details. 
\erema 


\section{Main Theorem}

\label{section:MainTheorem}

\theo \label{theorem:main1}
Let $\coneG, \ \coneK \subseteq \SymMat^n_+$ be closed convex cones such that 
$\coneG \subseteq \coneK$, $\Q \in \SymMat^n$, 
$\H \in \SymMat^n$, 
$\BC = \{\B^k : 1 \leq k \leq m\} \subseteq \SymMat^n$ and 
$\coneF_{\scriptsize \coneK}$ a face of $\coneK$ such that $\coneG\subseteq\coneF_{\scriptsize \coneK}$; hence  
$\coneG = \coneF_{\scriptsize \coneK}\cap\coneG$. 
Assume that conditions (C1), (C2), (C3), and (C4) below are satisfied. Then, 
COP$(\coneG\cap\coneJ_+(\BC),\Q,\H)$ is an exact convex relaxation of 
COP$(\bGamma^n\cap\coneG\cap\coneJ_+(\BC),\Q,\H)$.
\etheo
\noindent
\vspace{-8mm}
\begin{description} 
\item[(C1) ] If COP$(\coneG,\Q,\H)$ (equivalent to 
COP$(\coneF_{\scriptsize \coneK}\cap\coneG,\Q,\H)$) is solvable, then it is exact.  
\vspace{-2mm}
\item[(C2) ] 
$\coneJ_0(\B) \subseteq \coneK$ and $\coneF_{\scriptsize \coneK}\cap\coneJ_0(\B) \subseteq \coneG$ 
for every $\B \in \BC$. 
\vspace{-2mm}
\item[(C3) ] 
$\coneF_{\scriptsize \coneK}\cap\coneJ_0(\B) \subseteq \coneF_{\scriptsize \coneK} \cap\coneJ_+(\B')$
for every distinct $\B,\ \B' \in \BC$. 
\vspace{-2mm}
\item[(C4) ] COP$(\coneG\cap\coneJ_+(\BC),\Q,\H)$
(equivalent to COP$(\coneF_{\scriptsize \coneK}\cap\coneG\cap\coneJ_+(\BC),\Q,\H)$ 
has an optimal solution $\X$ at which 
{\em the KKT (Karush-Kuhn-Tucker) stationary condition} 
\begin{eqnarray}
\left. 
\begin{array}{l}
\X \in \coneF_{\scriptsize \coneK} \cap \coneG, \ \inprod{\B^k}{\X} \geq 0 \ (1 \leq k \leq m), \
\inprod{\H}{\X} = 1 \ \mbox{(primal feasibility)}, \\[3pt]
\displaystyle \bar{\y} \geq \0, \ \Q - \H \bar{t} -  \sum_{k=1}^m\bar{y}_k\B^k = \overline{\Y} \in 
(\coneF_{\scriptsize \coneK}\cap\coneG)^* 
\ \mbox{(dual feasibility)}, \\[3pt]
\displaystyle \sum_{k=1}^m\bar{y}_k \inprod{\B^k}{\X}= 0, \ 
\inprod{\overline{\Y}}{\X} = 0 \ \mbox{(complementarity)}
\end{array}
\right\}  
\label{eq:KKT1}
\end{eqnarray}
holds for some $(\bar{t},\bar{\y},\overline{\Y}) \in \Real \times \Real^m \times \SymMat^n$, 
where $(\bar{t},\bar{y},\overline{\Y})$ corresponds to an optimal solution of the 
dual of COP$(\coneF_{\scriptsize \coneK}\cap\coneG\cap\coneJ_+(\BC),\Q,\H)$.
\vspace{-1mm}
\end{description}

A proof of Theorem  \ref{theorem:main1} is given in Section~\ref{section:Proof}. 
The solvability and exactness of COP$(\coneK\cap\coneG,\Q,\H)$ in condition (C1) 
and COP$(\coneK\cap\coneG\cap\coneJ_+(\BC),\Q,\H)$ in (C4) 
can depend on $\Q \in \SymMat^n$ and 
$\H \in \SymMat^n$. In particular, we can take the exact convex relaxation of 
any QCQP listed in (a) through (g) in Section~\ref{section:Introduction} 
for COP$(\coneG,\Q,\H)$ in (C1). On the other hand, conditions (C2) and (C3) 
are independent from any choice of $\Q \in \SymMat^n$ and $\H \in \SymMat^n$. 
We note that every optimal solution $\X$ of COP$(\coneG\cap\coneJ_+(\BC),\Q,\H)$ 
satisfies~\eqref{eq:KKT1} under condition (C4). 

\medskip

Condition (C3) is similar to NIQC condition~\eqref{eq:condArimaMain1} with $\AC = \BC$. 
A difference is  that (C3) does not include `or $\coneF_{\scriptsize \coneK}\cap\coneJ_+(\B') \subseteq \coneF_{\scriptsize \coneK}\cap\coneJ_+(\B)$' 
for simplicity of the subsequent discussion; 
if it holds for distinct $\B, \ \B' \in \BC$, then we can delete $\B$ 
from $\BC$ in advance since 
$\coneG\cap\coneJ_+(\BC) = \coneG\cap\coneJ_+(\BC\backslash\{\B\})$. 
Another difference is that $\coneF_{\scriptsize \coneK}$ in (C3) is a face of a closed convex cone 
$\coneK \subseteq \SymMat^n_+$, while $\coneF$ in~\eqref{eq:condArimaMain1} is restricted to a face of $\SymMat^n_+$. 
We note that \eqref{eq:NIQCLMI0} with $\AC=\BC$ is a common sufficient condition 
for~\eqref{eq:condArimaMain1} 
and (C3) although their difference becomes critical  when 
{Theorem~\ref{theorem:main1} is applied} to exact CPP and DNN 
relaxations.
If $\coneG = \coneK = \coneF_{\scriptsize \coneK} = \SymMat^n_+$, then  
conditions (C1) and (C2) are satisfied and 
COP$(\coneJ_+(\BC),\Q,\H)$ becomes an exact SDP relaxation of 
COP$(\bGamma^n\cap\coneJ_+(\BC),\Q,\H)$ under conditions (C3) and (C4). 
This assertion corresponds to a special case of \cite[Theorem 1.5]{ARIMA2024}. 

\rema 
Assume that $\O \not= \H \in \SymMat^n_+$ and $(\coneK\cap\coneG\cap\coneJ_+(\BC))^* = (\coneK\cap\coneG)^* + \coneJ_+(\BC)^*$ 
(or equivalently $(\coneK\cap\coneG)^* + \coneJ_+(\BC)^*$ is closed). 
In this case, 
if COP$(\coneF\cap\coneG\cap\coneJ_+(\BC),\Q,\H)$ has an optimal solution $\X \in \coneK$, 
then~\eqref{eq:KKT1} holds.  See \cite[Theorem 1.1]{KIM2022}. 
\erema

\subsection{Some existing results related to Theorem~\ref{theorem:main1}}

\label{section:knownResults} 

Throughout this section, 
we assume that 
\begin{eqnarray*}
 & & \coneK = \SymMat^n_+, \  \BC : \mbox{a finite subset of $\SymMat^n$}, \
 \coneG = \coneJ_+(\AC) \ \mbox{for some finite subset $\AC$ of $\SymMat^n$}, \\[3PT]
& &  \coneF = \coneF_{\scriptsize \coneK} : \mbox{a face of $\SymMat^n_+$ that includes $\coneG = \coneJ_+(\AC)$}. 
  \end{eqnarray*}
 Then, both COP$(\coneG,\Q,\H)$ in  (C1) and COP$(\coneG\cap\coneJ_+(\BC),\Q,\H)$ 
  in (C4) can be written as COP$(\coneJ_+(\CC),\Q,\H)$ 
  \begin{eqnarray*}
  \eta(\coneJ_+(\CC),\Q,\H) & = & 
  \inf\left\{\inprod{\Q}{\X}: \X \in \coneJ_+(\CC), \ \inprod{\H}{\X} = 1\right\} \\
  & = & \inf\left\{\inprod{\Q}{\X}: \X \in \coneF\cap\coneJ_+(\CC), \ \inprod{\H}{\X} = 1\right\}. 
 \end{eqnarray*}
 Here either $\CC = \AC$ or $\CC = \AC\cup\BC$. 
 We also see that conditions (C2) and (C3) can be specialized and combined 
 to a single homogenized NIQC condition: 
 \begin{description}
  \item{(C2-3)' } $\coneF\cap\coneJ_0(\B) \subseteq \coneF\cap\coneJ_+(\B') \ \mbox{for every distinct } 
  \B \in \BC, \ \B' \in \AC\cup\BC$. 
 \end{description}
 
 
 Now, assume that $\H = \H^1$  
as in Section~\ref{section:standardQCQP}. Then, 
 COP$(\coneJ_+(\CC),\Q,\H^1)$ is an SDP relaxation of COP$(\bGamma^n\cap\coneJ_+(\CC),\Q,\H^1)$, represented as a QCQP: 
\begin{eqnarray}
  \eta(\bGamma^n\cap\coneJ_+(\CC),\Q,\H^1) & = & 
  \inf\left\{q(\u,\Q) : \u \in \Real^{n-1}, \ \u \in \CC_\geq \right\}  \label{eq:QCQPS2}  \\
  & = & \inf\left\{q(\u,\Q) : \u \in \Real^{n-1}, \ \u \in L_1(\coneF) \cap \CC_\geq \right\}. \nonumber 
 \end{eqnarray}
 
\lemm \label{lemma:SToR} \cite[Lemma 4.3]{ARIMA2024}
Let $\B,\ \B' \in \SymMat^n$. Assume that $\emptyset \not= L_1(\coneF)\cap\B_\leq\subseteq L_1(\coneF)\cap\B'_\geq$, then 
$\coneF\cap\coneJ_0(\B) \subseteq \coneF\cap\coneJ_-(\B) \subseteq \coneF\cap\coneJ_+(\B')$. 
\elemm
\noindent
The lemma above is slightly more general than the original \cite[Lemma 4.3]{ARIMA2024} 
where $\coneF = \SymMat^n_+$ is assumed. However, 
it can be derived from the original result since each face of $\SymMat^n_+$ is isomorphic to 
$\SymMat^r_+$ for some $r \in \{0,1,\ldots,n\}$, as mentioned in Remark~\ref{remark:facialReduction}.  
By Lemma~\ref{lemma:SToR}, conditions (C2-3)' can be replaced by a  non-homogenized NIQC condition: 
\begin{description}
\item{(C2-3)'' }  $\emptyset \not= L_1(\coneF)\cap\B_\leq \subseteq L_1(\coneF)\cap\B'_\geq  \ 
\mbox{for every distinct } \B \in \BC, \  \B' \in \AC \cup \BC$. 
 \end{description}
 It should be noted that condition (C2-3)'', which depends  on a special choice of $\H = \H^1$, is 
 sufficient for condition (C2-3)', which is independent of any choice of $\H \in \SymMat^n_+$. 
 
\medskip

 Let 
\begin{eqnarray*}
\spaceH^1 = \{ \X \in \SymMat^n : \inprod{\H^1}{\X} = 1 \} =  \{ \X \in \SymMat^n : X_{nn} = 1 \}. 
\end{eqnarray*}
Then, COP$(\bGamma^n\cap\coneJ_+(\CC),\Q,\H^1)$ and 
COP$(\coneJ_+(\CC),\Q,\H^1)$  can be rewritten as 
\begin{eqnarray*}
\eta(\bGamma^n\cap\coneJ_+(\CC),\Q,\H^1) & = & \inf 
\left\{ \inprod{\Q}{\X} : 
\X \in \bGamma^n\cap\coneJ_+(\CC)\cap \spaceH^1\right\}, 
\end{eqnarray*}
and 
\begin{eqnarray*}
\eta(\coneJ_+(\CC),\Q,\H^1) & = & \inf 
\left\{ \inprod{\Q}{\X} : \X \in \coneJ_+(\CC)\cap \spaceH^1\right\}, 
\end{eqnarray*}
respectively. 
Since the two problems above share the common linear objective function $\inprod{\Q}{\X}$ in 
$\X \in \SymMat^n$, the identity  
\begin{eqnarray}
& & \coneJ_+(\CC)\cap \spaceH^1 = 
\overline{\mbox{co}}(\bGamma^n\cap\coneJ_+(\CC)\cap \spaceH^1) \label{eq:convexHull}
\end{eqnarray}
on their feasible regions serves as a sufficient condition 
for the equivalence of 
COP$(\bGamma^n\cap\coneJ_+(\CC),\Q,\H^1)$ and COP$(\coneJ_+(\CC),\Q,\H^1)$, 
which played an essential role in \cite{JOYCE2024}. In particular,  
the following theorem mentioned in Section 1 was obtained from \cite[Corollaries 2]{JOYCE2024}  
as a special case, adapted to COP$(\bGamma^n\cap\coneJ_+(\CC),\Q,\H^1)$ and 
COP$(\coneJ_+(\CC),\Q,\H^1)$.  
\theo \label{theorem:JOYCE2024a} 
Assume that 
\begin{eqnarray}
& & 
\mbox{the identity \eqref{eq:convexHull} with $\CC = \AC$ holds}, \nonumber \\[3pt]
& &L_1(\coneF)\cap\B_=  \subseteq L_1(\coneF)\cap\B'_\geq \ \mbox{ for every distinct } \B \in \BC, \ \B' \in \AC\cup\BC, 
\label{eq:condOYCE2024a2}\\[3pt]
& & q(\cdot,\B) : \Real^{n-1} \rightarrow \Real \ \mbox{is not affine on $L_1(\coneF)$} \nonumber \\
& & \mbox{({\it i.e.}, the quadratic term of $q(\cdot,\B)|L_1(\coneF)$ is not identically zero)} \ (\B \in \BC).  
\label{eq:condOYCE2024a0}
\end{eqnarray}
\cite[Assumptions 3 and 4]{JOYCE2024}. 
Then the identity \eqref{eq:convexHull} with $\CC = \AC\cup\BC$ holds,
which implies $\eta(\coneJ_+(\AC\cup\BC),\Q,\H^1) =  \eta(\bGamma^n\cap\coneJ_+(\AC\cup\BC),\Q,\H^1)$ for every $\Q \in \SymMat^n$. \vspace{-1mm}
\etheo
\proof{
We provide a proof only for the case where $\coneF = \SymMat^n_+$ and $L_1(\coneF) = \Real^{n-1}$ 
since a general case can be reduced 
to this case by applying a linear isomorphism $\Phi$ from $\coneF$, 
which  includes $\coneJ_+(\AC)$,  onto $\SymMat^r_+$
for some $r \in \{0,1,\ldots,n\}$ as mentioned in Remark~\ref{remark:facialReduction}.
In this case, by \cite[Corollary 2]{JOYCE2024}, we have 
\begin{eqnarray*}
& &    \overline{\mbox{co}}(\bGamma^n\cap\coneJ_+(\AC\cup\BC)\cap\spaceH^1) =  \overline{\mbox{co}}(\bGamma^n\cap\coneJ_+(\AC)\cap\spaceH^1) \cap (\coneJ_+(\BC) \cap \spaceH^1)
\end{eqnarray*}
under assumption~\eqref{eq:condOYCE2024a2} and~\eqref{eq:condOYCE2024a0}. Hence it follows from 
the identity \eqref{eq:convexHull} with $\CC = \AC$ that 
\begin{eqnarray*}
& &    \overline{\mbox{co}}(\bGamma^n\cap\coneJ_+(\AC\cup\BC)\cap\spaceH^1) =  
(\coneJ_+(\AC) \cap \spaceH^1)\cap (\coneJ_+(\BC) \cap \spaceH^1) = (\coneJ_+(\AC \cup \BC) \cap \spaceH^1)
\end{eqnarray*}

\vspace{-7mm}

\begin{flushright} \qed \end{flushright}
}
%
\indent
Theorem~\ref{theorem:JOYCE2024a}  and Theorem~\ref{theorem:main1} 
with 
$\H = \H^1$ and $\coneG = \coneJ_+(\AC)$ 
are comparable. 
The main difference is: 
\mbox{the assumption \eqref{eq:convexHull} with $\CC = \AC$} in 
Theorem~\ref{theorem:JOYCE2024a} implies that 
$\eta(\coneJ_+(\AC),\Q,\H^1) =  \eta(\bGamma^n\cap\coneJ_+(\AC),\Q,\H^1)$ for 
any $\Q \in \SymMat^n$, while 
condition (C1) needs to be hold for a given $\Q \in \SymMat^n$. 
As mentioned in Section 1, this difference is critical in applications. 
Theorem~\ref{theorem:JOYCE2024a} applies only to classes of QCQPs whose SDP 
relaxations are exact independently  of the objective function,  
while  Theorem~\ref{theorem:main1} 
applies to QCQPs whose SDP relaxation exactness depends
 on both the objective and the constraint functions. 
Another difference 
lies in their non-homogenized NIQC conditions~\eqref{eq:condOYCE2024a2} and (C2-3)''. 
Clearly,  (C2-3)''  implies \eqref{eq:condOYCE2024a2}, however,  Theorem~\ref{theorem:JOYCE2024a} 
additionally imposes  condition~\eqref{eq:condOYCE2024a0}. In essence, conditions 
(C2-3)" and~\eqref{eq:condOYCE2024a2} with~\eqref{eq:condOYCE2024a0} are 
equivalent under reasonable assumptions, including Slater's constraint qualification 
and no-redundancy on $\AC\cup\BC$ in describing $\coneF\cap\coneJ_+(\AC\cup\BC)$. 
See \cite[Sections 4 and 5]{ARIMA2024} for more details.

\subsection{Proof of Theorem~\ref{theorem:main1}}
\label{section:Proof}

We need the following lemma for the proof. 

\lemm \label{lemma:YeZhang}
(\cite[Lemma 2.2]{YE2003}, see also \cite[Proposition 3]{STURM2003}) 
Let $\B \in \SymMat^n$ and $\overline{\X}\in \SymMat^n_+$ with 
rank$\overline{\X} = r \geq 1$. Suppose that $\inprod{\B}{\overline{\X}} \geq 0$. 
Then, there exists a rank-1 
decomposition of $\overline{\X}$ such that $\overline{\X} = \sum_{i=1}^r\x_i\x_i^T$ 
and $\inprod{\B}{\x_i\x_i^T} \geq 0$ 
$(1 \leq i \leq r)$. If, in particular, $\inprod{\B}{\overline{\X}}  = 0$, then 
$\inprod{\B}{\x_i\x_i^T} = 0$ $(1 \leq i \leq r)$. 
\elemm

\noindent
{\em Proof of Theorem~\ref{theorem:main1}}: 
Let $\X = \overline{\X}$ be an optimal solution of COP$(\coneF_{\scriptsize \coneK}\cap\coneG\cap\coneJ_+(\BC),\Q,\H)$ 
that satisfies the KKT condition~\eqref{eq:KKT1}, which serves as a sufficient condition for 
$\X$ to be an optimal solution of COP$(\coneF_{\scriptsize \coneK}\cap\coneG\cap\coneJ_+(\BC),\Q,\H)$. 
We have either case (i) $\inprod{\B^p}{\overline{\X}} = 0$ for some $p$ 
or case (ii) $\inprod{\B^k}{\overline{\X}} > 0$ $(1 \leq k \leq m)$.  
We first deal with case (i) $\inprod{\B^p}{\overline{\X}} = 0$.
Let $r = \mbox{rank}\overline{\X}$. By Lemma~\ref{lemma:YeZhang}, 
there exists a rank-1 decomposition of $\overline{\X}$ such that 
$\overline{\X} = \sum_{i=1}^r\x_i\x_i^T$ and 
$\x_i\x_i^T \in \coneJ_0(\B^p)$ $(1 \leq i \leq r)$. 
By (C2), $\x_i\x_i^T \in  \coneK \cap  \coneJ_0(\B^p) \subseteq \coneK \ (1 \leq i \leq r)$.
Since 
$\overline{\X} = \sum_{i=1}^r\x_i\x_i^T \in \coneF_{\scriptsize \coneK} \cap \coneG \cap \coneJ_0(\B) \subseteq 
\coneF_{\scriptsize \coneK}$ and $\coneF_{\scriptsize \coneK}$ is a face of $\coneK$, 
we see that $\x_i\x_i^T \in  \coneF_{\scriptsize \coneK}\cap\coneJ_0(\B^p)$. Hence   
$\x_i\x_i^T \in  \coneF_{\scriptsize \coneK}\cap\coneG\cap\coneJ_+(\BC)$ by (C2) and (C3) $(1 \leq i \leq r)$.   
Since 
$1 = \inprod{\H}{\overline{\X}} = \sum_{i=1}^r \inprod{\H}{\x_i\x_i^T}$, there exist a 
$\tau \geq 1/r$ and a $j \in \{1,\ldots,r\}$ such that $\inprod{\H}{\x_j\x_j^T}= \tau$. 
Let $\widetilde{\X} = \x_j\x_j^T/\tau$. 
Then $\widetilde{\X} \in  \coneF_{\scriptsize \coneK}\cap\coneG\cap\coneJ_+(\BC)$ and 
$\inprod{\H}{\widetilde{\X}} = \inprod{\H}{\x_j^T\x_j/\tau}= 1$. Hence $ \widetilde{\X}$ 
is a rank-$1$ feasible solution of COP$(\coneF_{\scriptsize \coneK}\cap\coneG\cap\coneJ_+(\BC),\Q,\H)$. Furthermore, 
we see from $\overline{\Y} \in (\coneF_{\scriptsize \coneK}\cap\coneG)^*$, $\bar{\y} \geq \0$ and 
$\widetilde{\X}, \  \x_i\x_i^T \in \coneF_{\scriptsize \coneK}\cap\coneG\cap\coneJ_+(\BC) \subseteq \coneF_{\scriptsize \coneK}\cap\coneG$ 
$(1 \leq i \leq r)$ that 
\begin{eqnarray*}
& & 
0 \leq \sum_{k=1}^m\bar{y}_k\inprod{\B^k}{\widetilde{\X}} = 
\sum_{k=1}^m\bar{y}_k \frac{\inprod{\B^k}{ \x_j\x_j^T}}{\tau} \leq
\sum_{k=1}^m\bar{y}_k \frac{\inprod{\B^k}{\displaystyle \sum_{i=1}^r\x_i\x_i^T}}{\tau} = 
\sum_{k=1}^m\bar{y}_k \frac{\inprod{\B^k}{\overline{\X}}}{\tau} = 0,  \\
& & 0 \leq \inprod{\overline{\Y}}{\widetilde{\X}} = \frac{\inprod{\overline{\Y}}{ \x_j\x_j^T}}{\tau} \leq 
\frac{\Inprod{\overline{\Y}}{\displaystyle \sum_{i=1}^r\x_i\x_i^T}}{\tau} = 
\frac{\inprod{\overline{\Y}}{\overline{\X}}}{\tau} = 0.  
\end{eqnarray*}
Therefore, $\X = \widetilde{\X}$ satisfies~\eqref{eq:KKT1}, and 
 is a rank-1 optimal solution of COP$(\coneF_{\scriptsize \coneK}\cap\coneG\cap\coneJ_+(\B),\Q,\H)$.   

\medskip

We now consider case (ii) $\inprod{\B^k}{\overline{\X}}  > 0$ $(1 \leq k \leq m)$. Then $\bar{\y}= 0$. 
 Hence 
 $\X  = \overline{\X}$ satisfies 
\begin{eqnarray*}
\begin{array}{l}
\X \in \coneF_{\scriptsize \coneK}\cap\coneG, \ \inprod{\H}{\X} = 1 \ \mbox{(primal feasibility)}, \\[3pt]
\displaystyle \Q - \H \bar{t}  = \overline{\Y} \in (\coneF_{\scriptsize \coneK}\cap\coneG)^*
\ \mbox{(dual feasibility)}, \\[3pt]
\inprod{\overline{\Y}}{\X} = 0 \ \mbox{(complementarity)}, 
\end{array}
\end{eqnarray*}
which serves as a sufficient condition for 
$\X \in \SymMat^n$ to be an optimal solution of COP$(\coneF_{\scriptsize \coneK}\cap\coneG,\Q,\H)$.  
Hence, $\overline{\X}$ is a common optimal solution of COP$(\coneF_{\scriptsize \coneK}\cap\coneG,\Q,\H)$ and 
COP$(\coneF_{\scriptsize \coneK}\cap\coneG\cap\coneJ_+(\BC),\Q,\H)$ with 
$\eta(\coneF_{\scriptsize \coneK}\cap\coneG,\Q,\H) = \eta(\coneF_{\scriptsize \coneK}\cap\coneG\cap\coneJ_+(\BC),\Q,\H) = \inprod{\Q}{\overline{\X}}$. 
By (C1), there exists a rank-$1$ optimal solution $\widehat{\X}$ of 
COP$(\coneF_{\scriptsize \coneK}\cap\coneG,\Q,\H)$, which satisfies 
\begin{eqnarray*}
	& & \inprod{\Q}{\widehat{\X}}= \inprod{\Q}{\overline{\X}},  \ 
	 \widehat{\X} \in \coneF_{\scriptsize \coneK}\cap\coneG, \ 
	\inprod{\H}{\widehat{\X}} = 1. 
\end{eqnarray*}
If $\inprod{\B^k}{ \widehat{\X}} \geq 0$ $(1 \leq k \leq m)$, then $\widehat{\X}$ is a rank-$1$ optimal 
solution of COP$(\coneF_{\scriptsize \coneK}\cap\coneG\cap\coneJ_+(\BC),\Q,\H)$. 
Otherwise, $\inprod{\B^q}{\widehat{\X}} < 0$ for some $q$. 
In this case, we can consistently define $\hat{\lambda} = \max \{ \lambda \in (0,1): 
\inprod{\B^k}{\lambda\widehat{\X} + (1-\lambda)\overline{\X}} \geq 0 \ (1 \leq k \leq m)\}$ since 
$\inprod{\B^k}{\overline{\X}} > 0$ $(1 \leq k \leq m)$ and $\inprod{\B^q}{\widehat{\X}} < 0$. 
Then $\widetilde{\X} = \hat{\lambda}\widehat{\X} + (1-\hat{\lambda})\overline{\X}$ is 
an optimal solution 
of COP$(\coneF_{\scriptsize \coneK}\cap\coneG\cap\coneJ_+(\BC),\Q,\H)$ such that $\inprod{\B^p}{\widetilde{\X}} = 0$ 
for some $p$. Thus, we have reduced this case to case (i).
\qed

\rema \label{numericalMethod} The proof of Lemma~\ref{lemma:YeZhang} given in 
\cite{YE2003} (see also \cite{STURM2003}) is constructive. Therefore, given a numerical method 
for solving COP$(\coneF_{\scriptsize \coneK}\cap\coneG,\Q,\H)$, the proof above provides a basic idea on how we compute 
a rank 1 optimal solution of COP$(\coneF_{\scriptsize \coneK}\cap\coneG\cap\coneJ_+(\BC),\Q,\H)$. 
\erema

\section{Exact SDP relaxations of QCQPs} 

\label{section:SDP} 

Throughout this section, we assume $\coneK = \SymMat^n_+$,  
and consider four different classes of QCQPs: 
QCQPs characterized by the rank-one generated 
property in Section~\ref{section:ROG}, 
convex QCQPs in Section~\ref{section:convex}, QCQPs characterized by sign pattern conditions 
in Section~\ref{section:signPattern} and continuous quadratic submodular minimization problems  
over the unit box 
in Section~\ref{section:submodular}. 
In each section, we present a simple example of 
$\AC, \ \BC \in \SymMat^3$, where conditions (C1) and (C2-3)" are satisfied 
with $n=3$, $\coneF = \SymMat^3_+$,  
$\coneG = \coneJ_+(\AC)$ 
and $\H = \H^1 \equiv \mbox{diag}(0,0,1)$.  
COP$(\coneJ_+(\AC),\Q,\H^1)$ (or COP$(\coneJ_+(\AC\cup\BC),\Q,\H^1)$) 
corresponds to the SDP relaxation of QCQP~\eqref{eq:QCQPS2} with $L_1(\coneF) = \Real^{2}$ 
and $\CC = \AC$ 
(or $\CC = \AC\cup\BC$, respectively). 
Thus, assuming that Condition (C4) holds with $\coneF_{\scriptsize \coneK} = \SymMat^3_+$, $\coneG = \coneJ_+(\AC)$,  
$\Q \in \SymMat^3$ and $\H = \H^1$, 
 COP$(\bGamma^n\cap\coneJ_+(\AC\cup\BC),\Q,\H^1)$ there 
can be solved by its SDP relaxation COP$(\coneJ_+(\AC\cup\BC),\Q,\H^1)$. 

\subsection{QCQPs characterized by the rank-one generated property}

\label{section:ROG}

The {\em rank-one-generated (ROG)} property \cite{ARGUE2023,ARIMA2024,KIM2020} 
is a well-known concept that provides a sufficient condition 
for the  equivalence of QCQPs and their convex relaxations. 
Let $\coneJ \subseteq \SymMat^n_+$ be a closed convex cone. 
$\coneJ$ is called ROG if it satisfies $\coneJ = \mbox{co}(\bGamma^n \cap \coneJ)$ 
or equivalently $\coneJ = \mbox{co}(\{\X \in \coneJ: \mbox{rank}\X \leq 1 \})$. 
Clearly, $\SymMat^n_+$ and $\coneCPP^n$ are ROG cones.
\lemm \label{lemma:ROG} \mbox{ \ } \vspace{-2mm}
\begin{description}
\item{(i) }  Let $\Q\in \SymMat^n$ and $\H \in \SymMat^n$. 
Assume that $\coneJ \subseteq \SymMat^n_+$ is a ROG cone. Then, 
\begin{eqnarray}
 -\infty < \eta(\coneJ,\Q,\H) \ \mbox{if and only if } -\infty < \eta(\coneJ,\Q,\H) =  \eta(\bGamma^n\cap\coneJ,\Q,\H), 
\label{eq:equiivalence}
\end{eqnarray}
(\cite[Theorem 3.1]{KIM2020},\cite[Lemma 20]{ARGUE2023}). \vspace{-2mm}
\item{(ii) } Every face of a ROG cone is ROG (\cite[Lemma 2.1 (iv)]{KIM2020},\cite[Lemma 3]{ARGUE2023}).
\end{description}
\elemm

\rema \label{remark:ROG} 
If a closed convex cone $\coneJ \subseteq \SymMat^n_+$ satisfies the assumption that 
$\O \not= \H \in \coneJ^*$ 
and $\Q -\H t$ lies in the interior of $\coneJ^*$ for some $t \in \Real$, then 
$-\infty < \eta(\coneJ,\Q,\H) < \infty$ implies the existence of an optimal solution $\X$ of 
COP$(\coneJ,\Q,\H)$ \cite[Theorem 2.1]{KIM2022}. If, in addition, $\coneJ$ is ROG, then 
COP$(\coneJ,\Q,\H)$ is exact 
({\it i.e.}, COP$(\coneJ,\Q,\H)$ has a rank-$1$ optimal solution) \cite[Theorem 2.4]{ARIMA2023}. 
For simplicity of discussion, we assume 
that if $-\infty < \eta(\coneJ,\Q,\H) =  \eta(\bGamma^n\cap\coneJ,\Q,\H) < \infty$ holds 
and COP$(\coneJ,\Q,\H)$ is solvable, 
then COP$(\coneJ,\Q,\H)$ is exact, although 
an additional assumption may be required.  
Thus, we assume that if $\coneG$ is ROG, then condition (C1) follows. 
\erema


Several sufficient conditions
on $\AC \subseteq \SymMat^n$ have been proposed in \cite{ARGUE2023,ARIMA2023,ARIMA2024} 
for $\coneJ_+(\AC)$ to be ROG.  
 Among them, our main theorem, Theorem~\ref{theorem:main1}, is motivated by the following results.  
 
 
 \theo 
 \label{theorem:ArimaMain1} 
 \cite[Theorems 1.5]{ARIMA2024}. 
 Let $\AC \subseteq \SymMat^n$ (possibly $|\AC| = \inf$) 
 and $\coneF$ be a face of $\SymMat^n_+$ 
 that includes $\coneJ_+(\AC)$. Assume that homogenized NIQC condition~\eqref{eq:condArimaMain1} holds. 
Then, $\coneJ_+(\AC) = \coneF\cap\coneJ_+(\AC)$ is ROG; hence \eqref{eq:equiivalence} holds with $\coneJ = \coneF \cap \coneJ_+(\AC)$. 
 \etheo
  

For COP$(\coneF\cap\coneG,\Q,\H)$ in condition (C1), we can take $\coneG = \coneJ_+(\AC')$ for some finite 
$\AC' \subseteq \SymMat^n$ and a face $\coneF$ of $\SymMat^n_+$ such that $\coneJ_+(\AC')\subseteq \coneF$ 
satisfying condition~\eqref{eq:condArimaMain1} with $\AC = \AC'$,  
and apply Theorem~\ref{theorem:main1}. In this case, the resulting 
$\coneG\cap\coneJ_+(\BC) = \coneJ_+(\AC'\cup\BC)$ satisfies condition~\eqref{eq:condArimaMain1} 
with $\AC = \AC' \cup \BC$. 
Hence, we could apply Theorem~\ref{theorem:ArimaMain1} directly to $\coneJ_+(\AC'\cup\BC)$.  
However, such an application of  Theorem~\ref{theorem:main1}  is not particularly interesting. Instead, 
we consider a condition on 
$\AC \subseteq \SymMat^n$, consisting of rank-2 matrices given in \cite[Section 3]{ARGUE2023} 
for $\coneJ_+(\AC)$ to be ROG, which is not covered by condition~\eqref{eq:condArimaMain1}.

\theo \label{theorem:Argue} \cite[Theorem 1]{ARGUE2023}. Let $\a \in \Real^n$,  
$\DC \subseteq \Real^n$. Assume that $\AC = \{ \a\d^T + \d\a^T : \d \in \DC\}$. Then, 
$\coneJ_+(\AC)$ is ROG; hence~\eqref{eq:equiivalence} holds with $\coneJ = \coneJ_+(\AC)$.  
\etheo 

\examp \label{example:Argue} 
Let $n=3$, $\coneF = \coneK = \SymMat^3_+$, $\Q \in \SymMat^3$, $\H^1 = \mbox{diag}(0,0,1)$, 
$\AC = \{\A^1, \A^2\}$, and $\BC = \{\B^1, \B^2\}$, where
\begin{eqnarray*}
& & 
\a = \begin{pmatrix} 1\\-2\\0\end{pmatrix}, \ 
\b = \begin{pmatrix}-2\\1\\0\end{pmatrix}, \
\cc = \begin{pmatrix} -1\\-1\\4\end{pmatrix}, \ 
\A^1= \a\b^T+\b\a^T = \begin{pmatrix}-4&5&0\\5&-4&0\\0&0&0\end{pmatrix}, \ \\[3pt]
& &  
 \A^2= \a\cc^T+\cc\a^T = \begin{pmatrix}-2&1&4\\1&4&-8\\4&-8&0\end{pmatrix}, \
 \B^1 = \begin{pmatrix}2&-2&1\\-2&2&1\\1&1&6\end{pmatrix}, \
\B^2 = \begin{pmatrix}1&0&-4\\0&1&-4\\-4&-4&31\end{pmatrix}. \
\end{eqnarray*} 
We consider COP$(\coneJ_+(\AC),\Q,\H^1)$, which is the SDP relaxation of 
QCQP~\eqref{eq:QCQPS2} with $\CC = \AC$, for condition (C1), and  
COP$(\coneJ_+(\AC\cup\BC),\Q,\H^1)$, which is the SDP relaxation of 
QCQP~\eqref{eq:QCQPS2} with $\CC = \AC\cup\BC$,  for condition (C4). 
By construction, $\AC$ satisfies the assumption of Theorem~\ref{theorem:Argue}. Hence, 
$\coneJ_+(\AC)$ is ROG, and (C1) is satisfied. 
We also see from Figure 2 that condition (C2-3)" is satisfied. 
\eexamp

\begin{figure}[t!]   \vspace{-35mm} 
\begin{center}
\includegraphics[height=110mm]{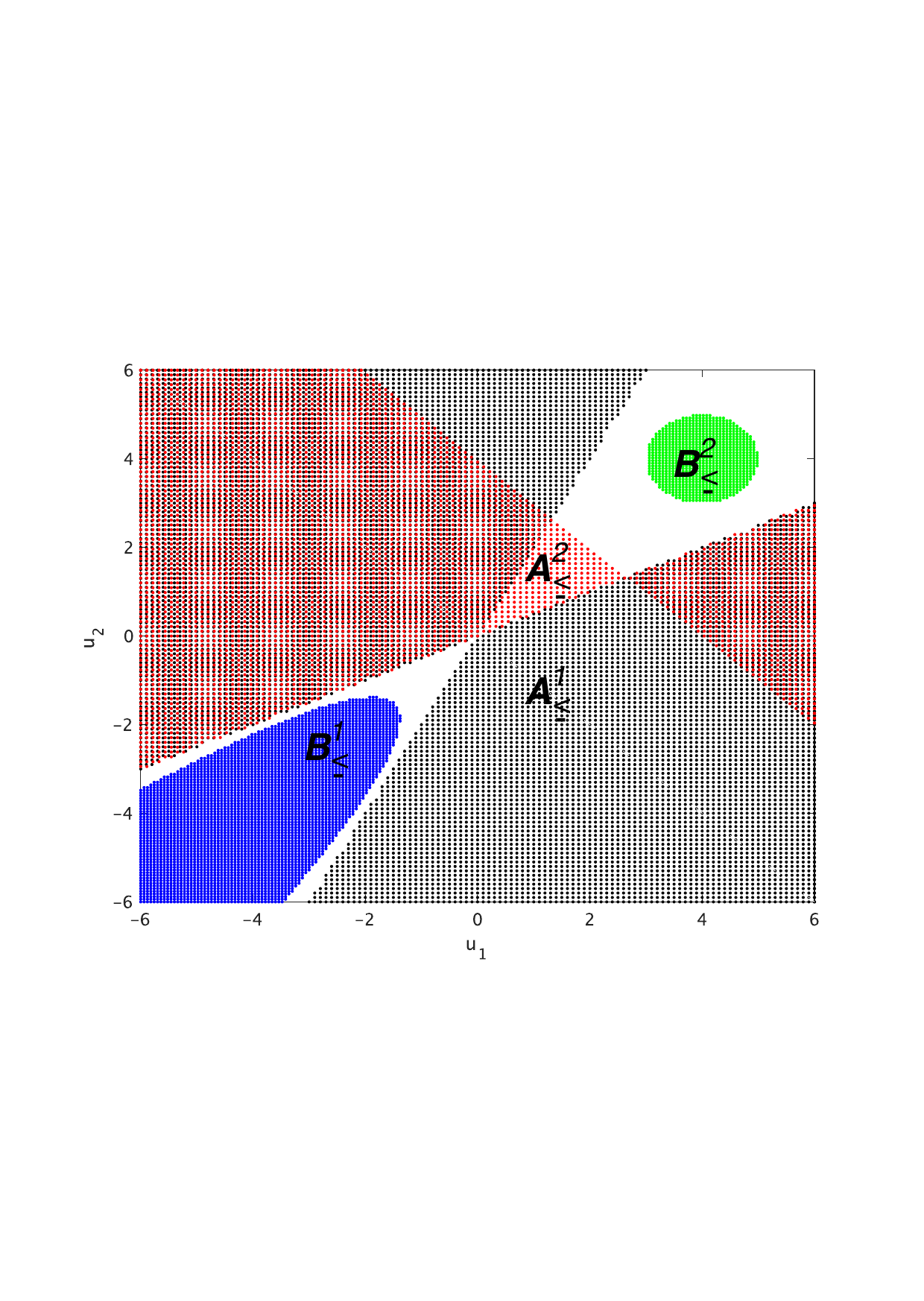}
\end{center}
\vspace{-30mm} 
\caption{
Example~\ref{example:Argue}. 
The unshaded (white) area corresponds to the interior of the feasible region $(\AC\cup\BC)_\geq$ of 
QCQP~\eqref{eq:QCQPS2} with $\CC = \AC\cup\BC$  and $L_1(\coneF) = \Real^{n-1}$. 
} 
\end{figure}


\subsection{Convex QCQPs} 

\label{section:convex}

The equivalence of a convex QCQP and its SDP relaxation 
 is well-known (\cite{SHOR1987},\cite[Section 4.2]{BENTAL2001}). 
Let $\AC$ be a finite subset of $\SymMat^n$ and $\Q \in \SymMat^n$.
We consider QCQP~\eqref{eq:QCQPS2} with $\CC = \AC$ and $L_1(\coneF) = \Real^{n-1}$, and assume 
that $q(\cdot,\Q) : \Real^{n-1} \rightarrow \Real$ and $-q(\cdot,\A):  \Real^{n-1} \rightarrow \Real$ 
$(\A \in \AC)$ are convex quadratic functions. Thus, \eqref{eq:QCQPS2} with $\CC = \AC$ 
 and $L_1(\coneF) = \Real^{n-1}$ forms a convex QCQP. 
 In this case, if its SDP relaxation, COP$(\coneJ_+(\AC),\Q,\H^1)$
is solvable, then it has rank-$1$ optimal solution. This fact is well-known and also easily proved. 

 \medskip
 
 \examp \label{example:convex}
Let 
\begin{eqnarray*}
& & 
\AC =\{\A^1,\A^2,\A^3\}, \ \BC = \{\B^1\},\ 
 \\[3pt]
& & \Q =  \begin{pmatrix} \C & \cc \\ \cc^T & 0 \end{pmatrix}\in \SymMat^3, \ \C \in \SymMat^2_+, \ 
\cc\in \Real^2, \ 
\A^1 = \begin{pmatrix}-1&0&0\\ 0&-1/2&0\\ 0&0&4 \end{pmatrix}, \\[3pt]
& & \A^2 = \begin{pmatrix}0&0&-1/2\\ 0&-1&0\\ -1/2&0&2 \end{pmatrix}, \
\A^3 = \begin{pmatrix}0&0&1/2\\ 0&-1&0\\ 1/2&0&4 \end{pmatrix}, \
 \B^1 = \begin{pmatrix}1/3&0&0\\ 0&1&0\\ 0&0&-1 \end{pmatrix}. \\
\end{eqnarray*}
We then see that $q(\cdot,\Q) : \Real^2 \rightarrow \Real, \ -q(\cdot,\A^k) : \Real^2 \rightarrow \Real$ 
$(1 \leq k \leq 3)$ are convex functions. Hence, condition (C1) is satisfied. 
We also see from Figure 3 that condition (C2-3)" is satisfied. 
\eexamp

\begin{figure}[t!]   \vspace{-30mm} 
\begin{center}
\includegraphics[height=120mm]{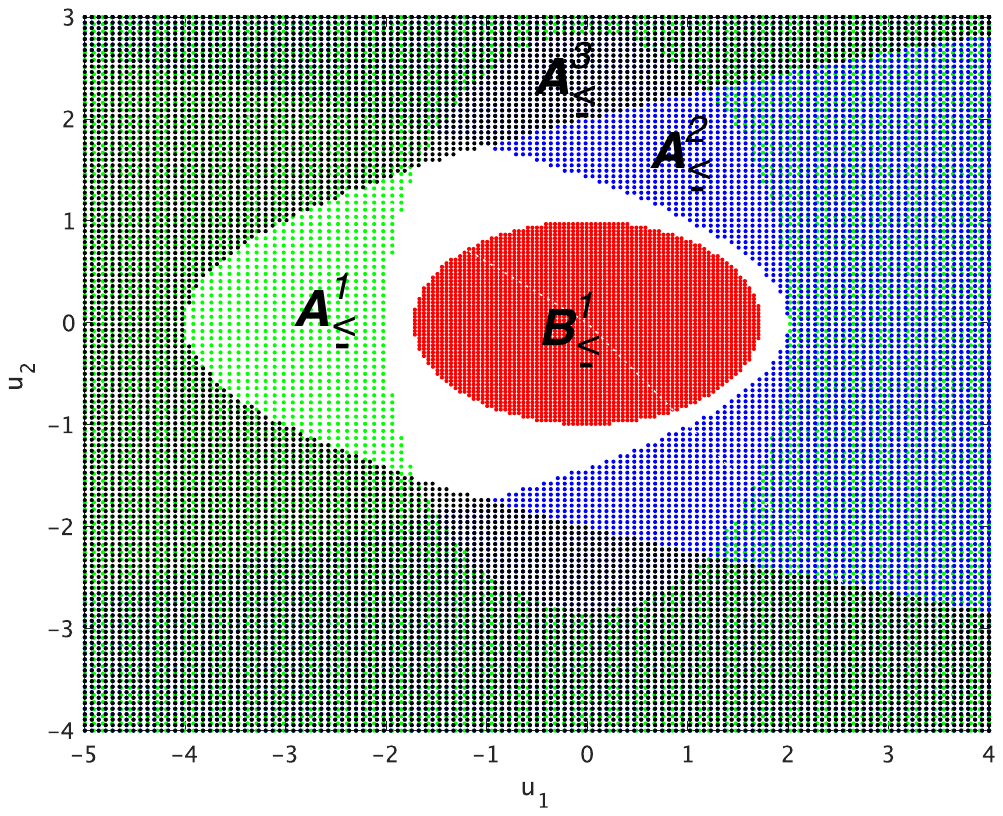}
\end{center}
\vspace{-35mm} 
\caption{Example~\ref{example:convex}. 
The unshaded (white) area represents the interior of the feasible region $(\AC\cup\BC)_\geq$ of 
QCQP~\eqref{eq:QCQPS2} with $\CC = \AC\cup\BC$  and $L_1(\coneF) = \Real^{n-1}$. 
} 
\end{figure}
 
\subsection{QCQPs characterized by sign pattern conditions 
 }

\label{section:signPattern}

QCQPs characterized by sign pattern conditions have been studied extensively in 
\cite{AZUMA2022,AZUMA2023,KIM2003,SOJOUDI2014}. 
In this section, we only consider 
a simple case to illustrate a  QCQP example satisfying conditions (C1) and (C2-3)". 
Specifically, we examine QCQP~\eqref{eq:QCQPS2} with 
$\CC = \AC = \{\A^k : 1 \leq k \leq \ell\}$. Let $\Q^0 = \Q, \ \Q^k = -\A^k \ (1 \leq k \leq \ell)$ and 
$L_1(\coneF) = \Real^{n-1}$. 
\lemm \label{lemma:signPattern} (\cite[Theorem 3.1]{KIM2003}, \cite[Corollary 1 (3)]{SOJOUDI2014}. 
Let $\H^1 = \mbox{diag}(0,\ldots,0,1) \in \SymMat^n_+$. Assume that $Q^k_{ij} \leq 0$ for every distinct $i, j \in \{1,\ldots,n\}$ $(0 \leq k \leq \ell)$. 
If COP$(\coneJ_+(\AC),\Q,\H^1)$ 
is solvable, then it has a rank-$1$ solution. 
\elemm 
\noindent 
See \cite[Theorem 2, Corollary 1]{SOJOUDI2014} for more general results. 

\examp \label{example:signPattern}  Let
\begin{eqnarray*}
& & 
n = 3, \ \AC =\{\A^1,\A^2,\A^3\} \subseteq \SymMat^3,\ \BC = \{\B^1,\B^2\} \subseteq \SymMat^3, \ 
\\[3pt]
& & \Q^0 \in \SymMat^3 \ \mbox{with all off-diagonal elements nonpositive, } 
\H = \mbox{diag}(0,0,1), \\
& & 
\Q^1 = -\A^1 = \begin{pmatrix}-1&0&0\\0&-1&-2\\0&-2&-3\end{pmatrix}, \
\Q^2 = -\A^2 = \begin{pmatrix}0&0&-1\\0&1&0\\-1&0&-6 \end{pmatrix}, \\
& & \Q^3 = -\A^3 = \begin{pmatrix}0&0&-2\\0&1&0\\-2&0&-4\end{pmatrix}, \
\B^1 = \begin{pmatrix}1 & 0& -3\\0&1&0\\-3&0&5\end{pmatrix}, \ 
\B^2  = \begin{pmatrix}0&0&-1\\0&2&0\\-1& 0& 10 \end{pmatrix}. 
\end{eqnarray*}
Obviously, all off-diagonal elements of $\Q^k$ $(0 \leq k \leq 3)$ are nonpositive. 
By Lemma~\ref{lemma:signPattern}, 
condition (C1) with $\coneF = \SymMat^3_+$ and $\coneG = \coneJ_+(\AC)$ is satisfied. 
We also see from Figure 4 that condition (C2-3)" is satisfied. 
\eexamp

\begin{figure}[t!]   \vspace{-40mm} 
\begin{center}
\includegraphics[height=150mm]{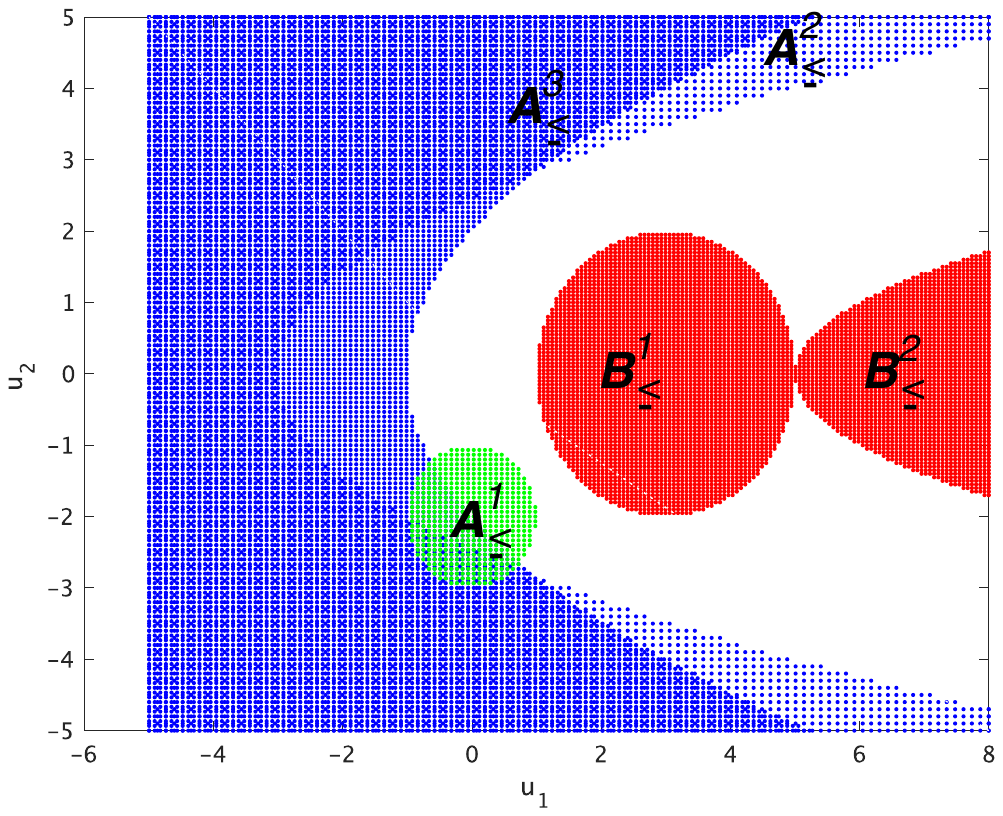}
\end{center}
\vspace{-45mm} 
\caption{
Example~\ref{example:signPattern}. 
The unshaded (white) area represents  the interior of the feasible region $(\AC\cup\BC)_\geq$ of 
QCQP~\eqref{eq:QCQPS2} with $\CC = \AC\cup\BC$  and $L_1(\coneF) = \Real^{n-1}$.  
} 
\end{figure}


\subsection{Continuous quadratic submodular minimization problem} 

\label{section:submodular}

Burer and Natarajan \cite{BURER2025} formulated 
the continuous quadratic submodular minimization problem over the box $[0,1]^{n-1}$ 
as
\begin{eqnarray*}
\varphi = \inf \left\{\u^T\C\u + \cc^T\u + \gamma : \u \in [0,1]^{n-1} \right\}, \label{eq:BurerQCQP}
\end{eqnarray*}
where $\C \in \SymMat^{n-1}$, $\cc \in \Real^{n-1}$ and $\gamma \in \Real$. 
As a convex relaxation of the problem, 
they proposed 
\begin{eqnarray}
\psi = \inf \left\{\inprod{\C}{\U} + 2\cc^T\u + \gamma : \U \leq \u\e^T, \ 
\begin{pmatrix} \U & \u \\ \u^T & 1 \end{pmatrix} \in \SymMat^n_+
 \right\}, \label{eq:BurerConvex}
\end{eqnarray}
where $\e$ denotes the $(n-1)$-dimensional column vector of $1$'s.

\theo \label{theorem:Burer} (\cite[Theorem 2]{BURER2025}). Assume that all off-diagonal elements 
of $\C \in \SymMat^{n-1}$ are nonpositive. Then, $\psi = \varphi$. 
\etheo 

We transform 
the convex relaxation~\eqref{eq:BurerConvex} into 
COP$(\coneJ_+(\AC),\Q,\H^1)$ for some $\AC \subseteq \SymMat^n$,  to which 
Theorem~\ref{theorem:main1} is then applied. 
Let $\Q = {\scriptsize \begin{pmatrix}\C&\cc\\ \cc^T & \gamma \end{pmatrix}}$ and 
$\H^1 = \mbox{diag}(0,\ldots,0,1)$.  
Define $\A^{pq} \in \SymMat^n$ $(1 \leq p \leq q \leq n-1)$ such that 
\begin{eqnarray*}
A^{pq}_{ij} = \left\{
\begin{array}{rl}
1/2 & \mbox{if $(i,j) = (p,n)$ or $(i,j) = (n,p)$}, \\
-1/2& \mbox{if $p\not=q$ and if $(i,j) = (p,q)$ or $(i,j) = (q,p)$}, \\
-1 & \mbox{if $p = q$, $i=p$ and $j = q$}, \\
0 & \mbox{otherwise}. 
\end{array}
\right. 
\end{eqnarray*}
Then, the constraint $\U \leq \u\e^T$ in~\eqref{eq:BurerConvex} can be written componentwisely as 
\begin{eqnarray*}
0 \leq u_p - U_{pq} = 
\inprod{\A^{pq}}{\X} \ (1 \leq p \leq q \leq n-1) \ \mbox{and } \inprod{\H^1}{\X} = 1,  
\end{eqnarray*}
where
$ \X = {\scriptsize \begin{pmatrix} \U & \u \\ \u^T & X_{nn} \end{pmatrix}}$. 
Hence, the convex relaxation~\eqref{eq:BurerConvex} is equivalent to COP$(\coneJ_+(\AC),\Q,\H^1)$, 
which satisfies Condition (C1) under the assumption of Theorem~\ref{theorem:Burer}. 
Also, we can easily verify that 
\begin{eqnarray*}
[0,1]^{n-1} \subseteq  \bigcap_{1 \leq p \leq q \leq n-1}  \left\{ \u \in \Real^{n-1} : 
\inprod{\A^{pq}}{{\scriptsize \begin{pmatrix} \u \\ 0 \end{pmatrix}}{\scriptsize \begin{pmatrix} \u \\ 0 \end{pmatrix}}^T} \geq 0 \right\} =  \bigcap_{1 \leq p \leq q\leq n-1} \A^{pq}_\geq  = \AC_\geq.  
\end{eqnarray*}

\examp \label{example:Burer}
Consider $\ell$ ellipsoids $\B^k_\leq$ 
 in $\Real^{n-1}$ represented by $\B^k \in \SymMat^n$ $(1 \leq k \leq \ell)$
such that 
\begin{eqnarray*}
& & \B^k_\leq = \left\{ \u \in \Real^{n-1} : \inprod{\B^k}{{\scriptsize \begin{pmatrix} \u \\ \1 \end{pmatrix}}{\scriptsize \begin{pmatrix} \u \\ \1 \end{pmatrix}}^T} \leq 0 \right\} \subseteq [0,1]^{n-1} \subseteq 
\AC_\geq, \\
& & \B^j_\leq \subseteq \B^k_\geq = 
\left\{ \u \in \Real^{n-1} : \inprod{\B^k}{{\scriptsize \begin{pmatrix} \u \\ 1 \end{pmatrix}}{\scriptsize \begin{pmatrix} \u \\ 1 \end{pmatrix}}^T} \geq 0 \right\} \ (j \not= k). 
\end{eqnarray*}
Let $\BC = \{\B^k : 1 \leq k \leq \ell\}$. Then, condition (C2-3)" is satisfied. 
\eexamp

\section{Exact CPP and DNN relaxations of QCQPs}

\label{section:CPPandDNN} 

In this section, we assume $\coneK \in \{\coneCPP^n,\coneDNN^n\}$, and discuss exact CPP and DNN 
relaxations of two classes of QCQPs.  The first  is a class of 
combinatorial QCQPs, presented in Section~\ref{section:binary}, and the second a class of standard quadratic optimization problems \cite{BOMZE2000,BOMZE2002},
discussed in Section~\ref{section:standardQOP}. 

\subsection{Combinatorial QCQPs} 

\label{section:binary}

Various 
combinatorial optimization problems,  including 
the maximum stable set problem, the quadratic unconstrained binary optimization problem (QUBO) 
and the quadratic assignment problem (QAP), are formulated as 
QCQPs with linear and complementarity constraints in nonnegative and binary variables. 
Burer \cite{BURER2009} demonstrated the exactness of the CPP relaxation for this class of QCQPs. 
Hence, this class of QCQPs can be used for condition (C1) with $\coneK = \coneCPP^n$.
However,  since $\coneCPP^n$ is numerically intractable,
the DNN relaxation \cite{GE2010,GOKMEN2022,KIM2013,KIM2016,YOSHISE2010}, though theoretically weaker than the CPP relaxation,
 is numerically tractable 
and often attains the exact optimal value of QCQP instances in this class. 
In fact, the Newton-bracketing method \cite{KIM2021}, based on DNN relaxation for solving QCQPs 
in the class, attains the optimal values of QUBO instances bqp100-2,$\ldots$,bqp100-5,bqp-8,bqp-9,bqp-10 
with dimension  100 \cite{WIEGELE2007}, and all 12 QAP instances chr with dimensions ranging from 12 to 
25 \cite{QAPLIB}. 
Therefore, Theorem~\ref{theorem:main1} can also be applied to such instances with $\coneK = \coneDNN^n$. 
We next provide additional details on the CPP and DNN relaxation of this class of QCQPs.

\medskip

Kim et al. \cite{KIM2020} showed that the class of QCQPs are characterized by the ROG property, 
and derived the exactness of their CPP relaxation from the ROG property. 
To illustrate the essential idea underlying their analysis,
we consider a simple problem of minimizing a quadratic function 
in three nonnegative variables $x_1,x_2,x_3$ subject to binary and complementarity constraints 
\begin{eqnarray*}
x_i \geq 0 \ (i=1,2,3), \ x_1(1-x_1) = 0 \ \mbox{and }x_2 x_3 = 0.  
\end{eqnarray*}
By introducing two additional nonnegative variables $x_4$ and $x_5$, 
this problem can be formulated as minimizing a quadratic form $\x^T\Q\x$ in $\x \in \Real^5$ for some 
$\Q \in \SymMat^5$ subject to the constraint 
\begin{eqnarray*}
& & x_i \geq 0 \ (i=1,2,3,4,5), \ x_1+x_4 - x_5 = 0, \ x_2x_3 = 0, \ x_1x_4 = 0. \ x_5^2 = 1, 
\label{eq:binary}
\end{eqnarray*}
or equivalently, 
$ \x  \in \Real^5_+, \  \x\x^T \in \coneJ_+(\AC) 
\ \mbox{and }   
\inprod{\H^1}{\x\x^T} = 1$, 
where 
\begin{eqnarray*}
& & \AC = \{\A^1,\A^2,\A^3\}, \ \H^1 = \mbox{diag}(0,0,0,0,1), \\ 
& & \A^1 = -{\scriptsize \begin{pmatrix}1\\0\\0\\1\\-1\end{pmatrix}}
{\scriptsize \begin{pmatrix}1\\0\\0\\1\\-1\end{pmatrix}}^T, \ 
\A^2 = -{\scriptsize \begin{pmatrix}0&0&0&0&0\\0&0&1&0&0\\0&1&0&0&0\\0&0&0&0&0\\0&0&0&0&0\end{pmatrix}},\ 
\A^3 = -{\scriptsize \begin{pmatrix}0&0&0&1&0\\0&0&0&0&0\\0&0&0&0&0\\1&0&0&0&0\\0&0&0&0&0\end{pmatrix}}.   
\end{eqnarray*}
As a result, we obtain 
QCQP~\eqref{eq:QCQPG} with $n=5$, $\coneK \in \{\coneDNN^5,\coneCPP^5\}$, 
$\coneG = \coneK\cap\coneJ_+(\AC)$  and $ \H = \H^1$, 
which is denoted as COP$(\bGamma^5\cap\coneK\cap\coneJ_+(\AC),\Q,\H^1)$ with 
$\coneK \in \{\coneDNN^5,\coneCPP^5\}$. 
COP$(\coneDNN^5\cap\coneJ_+(\AC),\Q,\H^1)$ and 
COP$(\coneCPP^5\cap\coneJ_+(\AC),\Q,\H^1)$ serve as its DNN relaxation and CPP relaxation, respectively. 
Since $-\A^k \in (\coneDNN^5)^* \subseteq (\coneCPP^5)^*$ $(k=1,2,3)$, 
$\coneDNN^5\cap\coneJ_+(\AC)$ and $\coneCPP^5\cap\coneJ_+(\AC)$ are faces 
of $\coneDNN^5$ and $\coneCPP^5$, respectively. 
The observation above can be generalized to formulate a higher-dimensional QCQP 
with linear and complementarity constraints in nonnegative and binary variables as 
COP$(\bGamma^n\cap\coneK\cap\coneJ_+(\AC),\Q,\H^1)$ 
with $\coneK \in \{\coneDNN^n,\coneCPP^n\}$ and some 
$\AC \subseteq \SymMat^n$ such that $\coneG = \coneK\cap\coneJ_+(\AC)$ forms 
a face of $\coneK$. 
See \cite{KIM2020} for more details. 
Therefore, QCQP and its convex relaxation under consideration are written as
COP$(\bGamma^n\cap\coneG,\Q,\H)$ and COP$(\coneG,\Q,\H)$, 
respectively, with $\H = \H^1$ and some face $\coneG$ of $\coneK \in \{\coneCPP^n,\coneDNN^n\}$. 

\medskip

Suppose that $\coneK = \coneCPP^n$. 
In  this case, 
$\coneCPP^n$ and its face $\coneG = \coneK\cap\coneJ_+(\AC)$ 
are ROG 
(Lemma~\ref{lemma:ROG} (ii)). 
By Lemma~\ref{lemma:ROG} (i), condition (C1) with $\coneF_{\scriptsize \coneK} = \coneG$ holds. 
Now, suppose that $\coneK = \coneDNN^n$.  
It is known that if $n \leq 4$, then $\coneCPP^n = \coneDNN^n$  
but if $n \geq 5$, then $\coneCPP^n$ is a proper subset of $\coneDNN^n$ \cite{BERMAN2003} 
and $\eta(\coneDNN^n\cap\coneJ_+(\AC),\Q,\H) = \eta(\coneCPP^n\cap\coneJ_+(\AC),\Q,\H)$ is 
not guaranteed. 
Even when $n \geq 5$, however, 
$\eta(\coneDNN^n\cap\coneJ_+(\AC),\Q,\H) = \eta(\coneCPP^n\cap\coneJ_+(\AC),\Q,\H)$ often 
holds for certain QCQP instances as  mentioned above.  For such instances, 
condition (C1) holds. 

\medskip

By letting $\coneF_{\scriptsize \coneK} = \coneG = \coneK\cap\coneJ_+(\AC)$, the second inclusion relation of 
condition (C2) is clearly satisfied. 
We next present an example of $\BC \subseteq \SymMat^n$ that satisfies the first inclusion 
relation of (C2) and (C3) for any face $\coneG$ of $\coneK \in \{\coneCPP^n,\coneDNN^n\}$. 

\begin{figure}[t!]   
\vspace{-30mm} 
\begin{center}
\includegraphics[height=120mm]{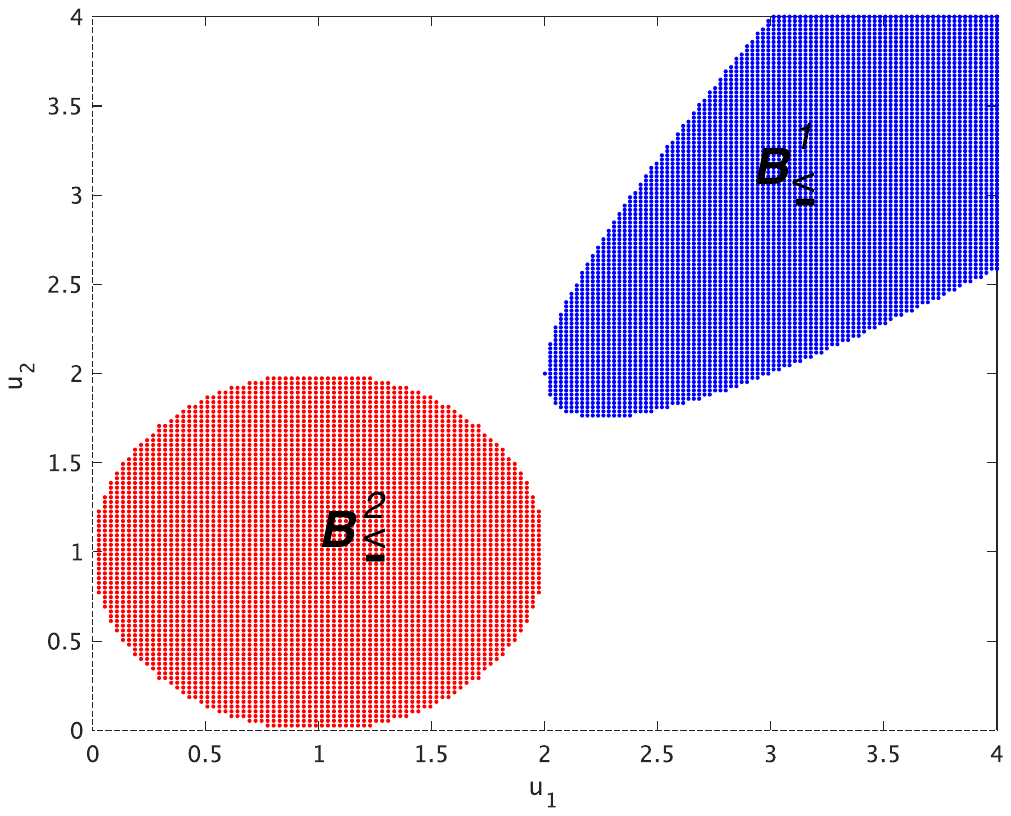}
\end{center}
\vspace{-35mm} 
\caption{
Example~\ref{example:CPP}: $n=3$, $\B^1_\leq = \{\u \in \Real^2 : -u_1+(u_2-u_1)^2 + 2 \leq 0\}$, 
$\B^2_\leq = \{ \u \in \Real^2 : \sum_{i=1}^{2} (u_i-1)^2 -1 \leq 0 \}$.} 
\end{figure}

\examp \label{example:CPP}

Let $\coneK\in \{\coneCPP^n,\coneDNN^n\}$, $\coneG = \coneF_{\scriptsize \coneK}$ a face of $\coneK$,  
$\Q \in \SymMat^n$ 
and $\H = \mbox{diag}(0,0,\ldots,1) \in \SymMat^n_+$. We consider the convex 
relaxation COP$(\coneG,\Q\,\H)$ of QCQP of the form~\eqref{eq:QCQPG}.
Condition (C1) is satisfied for any choice of $\Q \in \SymMat^n$ if $\coneK = \coneCPP^n$.  
We assume that (C1) is  also satisfied for $\coneK = \coneDNN^n$, as in the QUBO and QAP instances 
mentioned above. 
We provide a set $\BC = \{\B^1,\B^2\} \subseteq \SymMat^n$ that satisfies conditions (C2) and (C3). 
Choose  $\B^1, \ \B^2 \in \SymMat^n$ such that 
\begin{eqnarray*} 
\inprod{\B^1}{\x\x^T} 
& = & -x_1x_n +\sum_{i=2}^{n-1}(x_i-x_1)^2+2x_n^2, \\
\inprod{\B^2}{\x\x^T} 
& = &  \sum_{i=1}^{n-1}(x_i-x_n)^2  - x_n^2 
\end{eqnarray*}
for every $\x \in \Real^n$. See Figure 5.  It is easily verified that 
\begin{eqnarray*}
& & \B^1_= \subseteq \Real^{n-1}_+, \ 
\left\{\u  \in \Real^{n-1} : \inprod{\B^1}{{\scriptsize\begin{pmatrix}\u\\0\end{pmatrix}}
{\scriptsize\begin{pmatrix}\u\\0\end{pmatrix}}^T} = 0\right\}  \subseteq \Real^{n-1}_+ \cup (-\Real^{n-1}_+), \\
& & \B^2_=  \subseteq \Real^{n-1}_+, \ \left\{\u  \in \Real^{n-1} : \inprod{\B^2}{{\scriptsize\begin{pmatrix}\u\\0\end{pmatrix}}
{\scriptsize\begin{pmatrix}\u\\0\end{pmatrix}}^T} = 0\right\} 
= \{\0\}  \subseteq \Real^{n-1}_+ \cup (-\Real^{n-1}_+). 
\end{eqnarray*}
Hence, $\coneJ_0(\B^k) \subseteq \coneK$ $(k=1,2)$ holds by Lemma~\ref{lemma:C2} below, 
and the first inclusion relation of condition (C2) follows.
We also see that 
\begin{eqnarray*}
\emptyset \not= 
\B^1_\leq, \ \B^1_\leq \cap \B^2_< =  \emptyset, \ 
\emptyset \not= 
\B^2_\leq, \ \B^2_\leq \cap \B^1_< =  \emptyset, 
\end{eqnarray*}
which imply $\emptyset \not= \B^1_\leq \subseteq \B^2_\geq$ and $\emptyset \not= \B^2_\leq \subseteq \B^1_\geq$. 
By Lemma~\ref{lemma:SToR}, we obtain that $\coneJ_0(\B^1) \subseteq \coneJ_+(\B^2)$ 
and $\coneJ_0(\B^2) \subseteq \coneJ_+(\B^1)$. 
Hence condition (C3) is satisfied. 
Therefore, if condition (C4) holds, then 
the convex relaxation COP$(\coneG\cap\coneJ_+(\BC),\Q,\H^1)$ is exact  for 
COP$(\bGamma^n\cap\coneG\cap\coneJ_+(\BC),\Q,\H^1)$.
\eexamp

\lemm \label{lemma:C2} Let $\B \in \SymMat^n$. Then \vspace{-2mm}
\begin{eqnarray}
& & 
\coneJ_0(\B) \subseteq \coneDNN^n \Leftarrow \coneJ_0(\B) \subseteq \coneCPP^n \nonumber  \\
& &  \hspace{44mm} \Updownarrow \nonumber \\ 
& &  \hspace{22mm} \bGamma^n \cap \coneJ_0(\B) \subseteq \bGamma^n \cap \coneCPP^n 
\label{eq:NScondC1} \\
& & \hspace{44mm} \Updownarrow \nonumber \\
& & \hspace{14mm} \{ \x \in \Real^n:\x\x^T \in \coneJ_0(\B)\} \subseteq\Real^n_+\cup (-\Real^n_+) \label{eq:NScondC3} \\
& & \hspace{44mm} \Updownarrow \nonumber \\
& & \hspace{2mm} \B_= \subseteq \Real^{n-1}_+ \cup (-\Real^{n-1}_+) \ \mbox{and } \nonumber \\
& & \hspace{20mm} \left\{\u \in \Real^{n-1} : \inprod{\B}{
{\scriptsize \begin{pmatrix}\u\\0\end{pmatrix}}
{\scriptsize \begin{pmatrix}\u\\0\end{pmatrix}}^T} = 0 \right\} 
\subseteq \Real^{n-1}_+ \cup (-\Real^{n-1}_+). 
\label{eq:NScondC4}
\end{eqnarray}
\elemm
\proof{The first $\Leftarrow$ is straightforwrd since $\coneCPP^n \subseteq \coneDNN^n$. 
$\coneJ_0(\B) \subseteq \coneCPP^n \Rightarrow $ \eqref{eq:NScondC1} is also obvious. 
To prove the converse, we apply Theorem~\ref{theorem:ArimaMain1} by 
taking $\AC= \{\B,-\B\}$ and $\coneF = \SymMat^n_+$ to see that $\coneJ_0(\B) = \coneJ_+(\AC)$ 
is ROG. Also $\coneCPP^n$ is ROG. Hence, taking the convex hull of 
$\bGamma^n \cap \coneJ_0(\B) \subseteq \bGamma^n \cap \coneCPP^n$, $\coneJ_0(\B) \subseteq \coneCPP^n$ follows. \eqref{eq:NScondC1} $\Leftrightarrow$ \eqref{eq:NScondC3} and \eqref{eq:NScondC3} $\Rightarrow$ \eqref{eq:NScondC4} are easily verified. To prove \eqref{eq:NScondC3} $\Leftarrow$ \eqref{eq:NScondC4}, 
assume that \eqref{eq:NScondC4} holds. Let $\x = {\scriptsize \begin{pmatrix} \u\\z\end{pmatrix}} \in 
\{ \x \in \Real^n:\x\x^T \in \coneJ_0(\B)\}$ or equivalently $\inprod{\B}{
{\scriptsize \begin{pmatrix} \u\\z\end{pmatrix}}{\scriptsize \begin{pmatrix} \u\\z\end{pmatrix}}^T}
= 0$. If $z = 0$ then ${\scriptsize \begin{pmatrix} \u\\z\end{pmatrix}} \subseteq \Real^n_+\cup (-\Real^n_+)$ 
follows from the latter inclusion relation of \eqref{eq:NScondC4}. Otherwise 
$\inprod{\B}{
{\scriptsize \begin{pmatrix} \u/z \\1\end{pmatrix}}{\scriptsize \begin{pmatrix} \u/z\\1\end{pmatrix}}^T}
= 0$. Hence, $\u/z \in \B_= \subseteq \Real^{n-1}_+ \cup (-\Real^{n-1}_+)$, which implies 
${\scriptsize \begin{pmatrix} \u\\z\end{pmatrix}} \in \Real^n_+ \cup (-\Real^n_+)$. 
\qed
}

\subsection{The standard quadratic optimization problem} 

\label{section:standardQOP}

Let $\Q \in \SymMat^n$. We consider the problem of minimizing the quadratic form $\inprod{\Q}{\x\x^T}$ 
over the standard simplex:
\begin{eqnarray*}
\varphi & = & \inf\left\{ \inprod{\Q}{\x\x^T} : \x \in \Real^n_+, \ \e^T\x = 1 \right\}, 
\end{eqnarray*}
which is called {\em the standard quadratic optimization problem}. Here, $\e$ denotes the $n$-dimensional 
vector of $1$'s.  It is well-known that the optimal value $\varphi$ is nonnegative if and only if $\Q$ is copositive \cite{DUR2021}.  
By letting $\H = \e\e^T$ and $\coneG = \coneK \in \{\coneCPP^n,\coneDNN^n\}$, 
we represent the problem as COP$(\bGamma^n\cap\coneK,\Q,\H)$ of the form 
QCQP~\eqref{eq:QCQPG} with $\coneG = \coneK$. 
In the case where $\coneK = \coneCPP^n$, we have the exact CPP relaxation 
COP$(\coneCPP^n,\Q,\H)$, 
which satisfies condition (C1) with $\coneG = \coneF_{\scriptsize \coneK} = \coneCPP^n$, $\H = \e\e^T$ and 
every $\Q \in \SymMat^n$   
\cite{BOMZE2000}. 

\medskip

In the case where $\coneK = \coneDNN^n$, we have the DNN relaxation COP$(\coneDNN^n,\Q,\H)$, which is 
not necessarily exact.  
In the recent paper \cite{GOKMEN2022}, G\"{o}kmen and Yildirim provided 
some sufficient conditions on $\Q$ for 
$\eta(\coneDNN^n,\Q,\H) = \eta(\coneCPP^n,\Q,\H)  =\varphi$, which 
include conditions
\begin{eqnarray*}
& & \Q \in \QC^1 \equiv \left\{ \Q \in \SymMat^n :  
\displaystyle \mbox{min}\{Q_{ij}: 1\leq i,j \leq n\} = \mbox{min}\{ Q_{kk} : 1\leq k \leq n  \} \right\}, \\
& & \Q \in \QC^{\rm concave} \equiv \left\{ \Q \in \SymMat^n :  \inprod{\Q}{\d\d^T} \leq 0 \ \mbox{if } \e^T\d = 0 \right\}, \\ 
& & \Q \in \QC^{\rm convex} \equiv \left\{ \Q \in \SymMat^n :  \inprod{\Q}{\d\d^T} \geq 0 \ \mbox{if } \e^T\d = 0 \right\}. 
\end{eqnarray*}
We note that $\QC^{\rm concave} \subseteq \QC^1$ was shown there. They also stated another condition 
induced from maximum weighted cliques in perfect graphs (see Section 4.3 of \cite{GOKMEN2022}). 
If $Q \in \SymMat^n$ satisfies one of those conditions, 
COP$(\coneDNN^n,\Q,\H)$ is exact; 
hence condition (C1) holds with $\coneG = \coneF_{\scriptsize \coneK}  = \coneDNN^n$ 
and $\H = \e\e^T$. 

\medskip

In both cases $\coneK=\coneCPP^n$ and $\coneK=\coneDNN^n$ above, $\BC \subseteq \SymMat^n$ 
given in Example~\ref{example:CPP} satisfies 
conditions (C2) and (C3) with $\coneG = \coneF_{\scriptsize \coneK} = \coneK$.

\section{Conclusions}

\label{section:Conclusion}

We have shown how a QCQP, whose convex relaxation is known to be exact, can be systematically expanded by 
adding non-intersecting quadratic inequalities. 
In particular, we have established Theorem~\ref{theorem:main1}, which provides a set of sufficient conditions (C1), (C2), (C3), and (C4), 
ensuring  that  the convex relaxation of the expanded QCQP  remains exact. Condition (C1) requires the convex relaxation of the original QCQP 
to be exact if it is solvable. 
Conditions (C2) and (C3) are homogenized NIQC conditions on the added quadratic constraints. 
These two conditions may be replaced by non-homogenized NIQC condition (C2-3)" 
for SDP relaxations of the standard-form QCQP~\eqref{eq:QCQPS}.
Condition (C4) 
requires that the convex relaxation of the expanded QCQP has an optimal solution that 
satisfies the KKT stationary condition. Among these conditions, (C1) is necessary and (C4) 
is natural,
whereas the homogenized NIQC conditions (C2) and (C3) are essential for maintaining the exactness in the expanded formulation. 
Consequently,  homogenized NIQC conditions (C2) and (C3) play  a central  role 
in Theorem~\ref{theorem:main1}. 


\end{document}